\documentclass{amsart}
\usepackage[utf8]{inputenc}
\usepackage[english]{babel}

\usepackage{amsthm}
\usepackage{amsmath}
\usepackage{amssymb}

\usepackage{graphicx}
\usepackage{commath}
\usepackage{breqn}
\usepackage{xcolor}
\usepackage{comment}
\usepackage{hyperref}

\usepackage{tikz}
\usetikzlibrary{decorations.pathmorphing}
\usetikzlibrary{decorations.pathreplacing}
\usetikzlibrary{arrows}
\usetikzlibrary{arrows, arrows.meta}
\usepackage{caption}
\tikzset{snake it/.style={decorate, decoration=snake}}

\usepackage{lmodern}

\DeclareMathOperator{\diam}{diam}
\DeclareMathOperator{\Cay}{Cay}

\title{Random divergence of groups}

\usepackage{dsfont}

\newtheorem{defn}{Definition}[section]
\newtheorem{lemma}[defn]{Lemma}
\newtheorem{thm}[defn]{Theorem}
\newtheorem{claim}{Claim}
\newtheorem{prop}[defn]{Proposition}
\newtheorem{notation}{Notation}

\newtheorem{Assumpt}{Assumption}
\newtheorem{cor}[defn]{Corollary}
\newtheorem{example}[defn]{Example}
\newtheorem{remark}{Remark}
\newtheorem{question}[]{Question}
\newcommand{\dive}{{\mathrm{div}}}

\newcommand{\nest}{\sqsubseteq}

\newcommand{\tsh}[1]{\left\{\kern-.7ex\left\{#1\right\}\kern-.7ex\right\}}
\newcommand{\Tsh}[2]{\tsh{#2}_{#1}}
\newcommand{\ignore}[2]{\Tsh{#2}{#1}}

\newtheorem{assumptx}{Assumption}

\usepackage{bbold}
\def\acts{\curvearrowright}

\usepackage{amssymb}
\usepackage{mathtools} 
\usepackage{float}

\theoremstyle{plain}

	\author[Antoine Goldsborough]{Antoine Goldsborough}
	\address{Maxwell Institute and Department of Mathematics, Heriot-Watt University, Edinburgh, UK}
	\email{ag2017@hw.ac.uk}

\author[Alessandro Sisto]{Alessandro Sisto}
	\address{Maxwell Institute and Department of Mathematics, Heriot-Watt University, Edinburgh, UK}
	\email{a.sisto@hw.ac.uk}

\begin{document}

\maketitle

\begin{abstract}
The divergence of a group is a quasi-isometry invariant defined in terms of pairs of points and lengths of paths avoiding a suitable ball around the identity. In this paper we study ``random divergence'', meaning the divergence at two points chosen according to independent random walks or Markov chains; the Markov chains version can be turned into a quasi-isometry invariant. We show that in many cases, such as for relatively hyperbolic groups, mapping class groups, and right-angled Artin groups, the divergence at two randomly chosen points is with high probability equivalent to the divergence of the group. That is, generic points realise the largest possible divergence.
\end{abstract}

\section{Introduction}
The notion of divergence of groups was introduced and studied in \cite{Gromov93} and \cite{Gersten}, and roughly it measures the lengths of paths joining two points $x,y$, and avoiding a ball of suitable radius around the identity, as a function of the distance of $\{x,y\}$ from the identity. It arose in the study of non-positively curved manifolds and metric spaces. As it turns out, up to equivalence of functions, the divergence of a group is a quasi-isometry invariant. 

The divergence of a group has been computed in many examples, we recall some here. In \cite{GerstenDivergence}, Gersten showed that the divergence of graph manifold groups is quadratic, and in fact this characterises graph manifolds among closed $3$-manifolds. 
Mapping class groups also have quadratic divergence except for a few exceptional surfaces \cite{BehrsotckMCG,DuchinRafi}.
It was shown in \cite{Macura, BD:thick} that for any degree $d$, there exists a CAT(0) group whose divergence is a polynomial of degree $d$, and a similar result was also obtained in the case of right-angled Coxeter groups \cite{DaniThomas}.
It was shown in \cite{SistoRelativeHyp} that finitely presented one-ended relatively hyperbolic groups have exponential divergence. Finally, groups with exotic behaviour of the divergence function are constructed in \cite{OOS} and \cite{GruberSisto}.

The divergence measures the ``worst-case scenario'' for points $x,y$ as above, meaning the points that are ``hardest'' connect outside a ball around the identity (given a bound on the distance from the identity). It is natural to wonder whether two ``generic'' points are easier to connect than the worst case scenario or not. This can be made precise by choosing the points using independent random walks. In fact, using the point of view from \cite{GS21} of Markov chains as a ``quasi-isometry invariant'' version of random walks, it is natural to consider Markov chains as well in this context. Indeed, considering divergence of pairs of points chosen using independent Markov chains yields quasi-isometry invariants in a few possible ways (e.g. taking expectactions). We do not pursue this here as actually in this paper we show that in many cases generic points realise the worst case scenario for divergence, as we now state. We denote $\dive(x,y,p;\delta)$ the divergence of points $x,y$ with respect to the basepoint $p$ and auxiliary constant $\delta$, see Definition \ref{defn:divergence} for the precise definition. Also, the divergence of a group is only defined up to a certain equivalence relation on functions that we denote $\asymp$.

\begin{thm}
\label{thm:intro}
    Let $(w_n)_n,(z_n)_n$ be independent copies of a random walk driven by a measure whose finite support generates $G$ as a semigroup, where $G$ is one of the following:
\begin{itemize}
    \item a finitely presented one-ended relatively hyperbolic group with infinite index peripheral subgroups.
    \item the fundamental group of a non-geometric graph manifold.
    \item the mapping class group of a closed connected oriented surface of genus at least 2.
    \item a right-angled Artin group whose defining graph is connected and not a join.
\end{itemize}
 Then for all sufficiently small $\delta>0$ there exists a function $g$ with $g(n)\asymp \dive_G(n)$ such that
$$\mathbb P[\dive(w_n,z_n,1;\delta)\geq g(n)]\to 1.$$
Moreover, the same is true in the first two cases replacing the random walks with tame Markov chains in the sense of \cite{GS21}.
\end{thm}

We note that work in progress of the authors with Mark Hagen and Harry Petyt will yield a Markov chain version of the statement for the remaining two cases as well.

The theorem follows from our main theorem, Theorem \ref{thm:lower_bound}, which roughly speaking gives a lower bound on random divergence for many more cases than the ones mentioned above as well, and above we chose some significant examples where this lower bound matches the divergence of the group. One example that is covered by Theorem \ref{thm:lower_bound} is that of acylindrically hyperbolic groups (in view of Lemma \ref{lem:fdivelement_acylindrically_hyperbolic}), yielding the following.

\begin{thm}
    Let $G$ be an acylindrically hyperbolic group and let $(w_n)_n,(z_n)_n$ be independent copies of a random walk driven by a measure whose finite support generates $G$ as a semigroup. Then for all sufficiently small $\delta>0$ there exists a superlinear function $g$ such that
$$\mathbb P[\dive(w_n,z_n,1;\delta)\geq g(n)]\to 1.$$
\end{thm}

While we find it interesting that generic points have the largest possible divergence in many cases, this is in a sense a negative result, meaning that this point of view does not yield new quasi-isometry invariants in those cases. It is natural to ask:

\begin{question}
Does there exist a group which has random divergence (in any of its possible meanings) which is strictly lower than its divergence?
\end{question}

We believe that both a positive and a negative answer would be interesting. If the answer is positive, then random divergence is a new quasi-isometry invariant. If the answer is negative, then in \emph{every group} generic points realise the largest divergence, which would be very surprising.

\subsection*{Outline}

In Section \ref{sec:f_div} we introduce the notion of $f$-divergent element of a group, where $f$ is a function, which is similar to a notion introduced in \cite{GS21} and has implicitly been used in certain computations of divergence, for example for mapping class groups. Roughly, $f$-divergence gives lower bounds on length of paths outside neighborhoods of the cyclic group generated by the element, in terms of the function evaluated at the radius of the neighborhood. Lemma \ref{lem:lower_bound_divergence} says that if a group contains an $f$-divergent element then the divergence of the group is $\succeq nf(n)$, and this bound is optimal in the cases mentioned in Theorem \ref{thm:intro}. We note that any acylindrically hyperbolic group contains an $f$-divergent element for some diverging function $f$.

In Section \ref{sec:statements} we state the various assumptions under which we can control random divergence, and in fact we have two different sets of assumptions for random walks and for Markov chains, Assumptions \ref{assump:random_walk} and \ref{assump:markov_chain}. In this section we also state our main theorem, Theorem \ref{thm:lower_bound}. In fact, we decided to turn two key propositions in the proof of the theorem into black boxes, so that in the future to extend the theorem to other cases one only needs to prove the two propositions. The scheme of the proof is the following:
\begin{figure}[ht]
    \centering
    \includegraphics[scale=0.3]{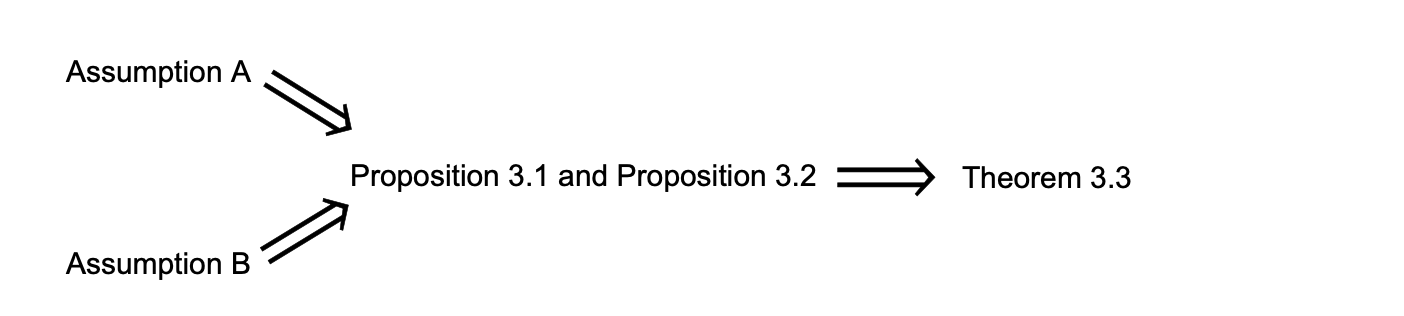}
    \caption{A recap of implications in this paper }
    \label{fig:my_label}
\end{figure}

In Section \ref{sec:probability} we establish the propositions under the two sets of assumptions, and in Section \ref{sec:main_proof} we prove Theorem \ref{thm:lower_bound}. Using this and computations of divergence from the literature we then prove Theorem \ref{thm:intro}.

\subsection*{Acknowledgments} This work was started at the AIM Workshop ``Random walks beyond hyperbolic groups'' organised by Joseph Maher, Yulan Qing, and Giulio Tiozzo. We thank the organisers and AIM, as well as all the participants who discussed various aspects of random divergence at the workshop and beyond. We thank in particular Michael Hull and Jing Tao for working on random divergence with us since the workshop.

\section{$f$-divergence}
\label{sec:f_div}

In this section, we define $f$-divergent elements and give some examples of groups containing such elements. The motivation for this study is that often, the divergence of a group originates from the existence of such an element.

\subsection{Definition and examples}

We first define an $X$-projection; roughly the idea is that given a group $G$ acting on a hyperbolic space $X$, we can use closest point projections in $X$ to define coarse retractions onto subsets of $G$.

\begin{defn}(\cite[Definition 3.2]{GS21})
\label{defn:X_proj}
	Let $G$ act on a hyperbolic space $X$, with a fixed basepoint $x_0$, and let $A \subseteq G$. An $X$-projection $\pi:G \to A$ is a map such that for all $h \in G$, we have that $\pi(h)x_0$ in a closest point in $Ax_0$ to $hx_0$. 
\end{defn}

We can now give the definition of an $f$-divergent element.

\begin{figure*}[h!]
	\centering
	\begin{tikzpicture}[scale=0.8]
\draw[black] (0,0) to[out=20,in=160] node [pos=0.75] (mid1) {} 
        node [pos=0.55] (mid2) {} 
         node [pos=0.85] (mid3) {} (10,0);
\filldraw[black] (2,0.63) circle (1pt) node[anchor=north]{$\pi(x)$};    
\filldraw[black] (8,0.63) circle (1pt) node[anchor=north]{$\pi(y)$}; 

\draw[densely dashed][red] (0,1.5) to[out=20,in=160] node [pos=0.75] (mid1) {} 
        node [pos=0.55] (mid2) {} 
         node [pos=0.85] (mid3) {} (10,1.5);     

\filldraw[black] (1,3) circle (1pt) node[anchor=east]{$x$};
\filldraw[black] (9,3) circle (1pt) node[anchor=west]{$y$};

\draw[-stealth] (1,3) --(1.97,0.68);
\draw[-stealth] (9,3) --(8,0.68);

\draw [->,decorate,decoration=snake] (1,3) -- (9,3);
  \path [draw=blue,snake it] (1,3) -- (9,3);

\draw[blue] (6,3.4) node{$\ell(\alpha) >f(d)$};

\draw[>=triangle 45, <->,purple] (2.3,0.6)-- (7.7,0.6) ;

\draw[purple] (5,0.3) node{$\geq \theta$};

\draw[red] (10,1.2) node{$\mathcal N_{d}(A)$};
\draw[black] (10,-0.3) node{$A$};

\end{tikzpicture}
\caption{Picture of $f$-divergence taken from \cite{GS21} (where $f$ was assumed to be super-linear, but this is not the case here).}\label{fig:f_divergent}
\end{figure*}
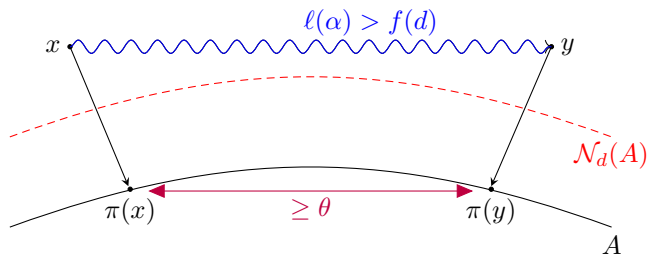

\begin{defn}
\label{defn:f-divergent}
Let $G$ be a group with a fixed word metric.
Let $f:\mathbb R^+\to \mathbb R^+$ be an increasing diverging function. We say that $g\in G$ is $f$-divergent if there exists a non-elementary action of $G$ on some hyperbolic space $X$ where $g$ is loxodromic, and a constant $\theta$ with the following property. Let $\pi:G\to \langle g \rangle$ be a fixed $X$-projection. If $x,y\in G$ have $d_G(\pi(x),\pi(y))\geq \theta$ and $\alpha$ is a path in $G$ from $x$ to $y$ avoiding the $d$-neighborhood of $\langle g \rangle$, then $\alpha$ has length at least $f(d)$. See Figure \ref{fig:f_divergent}.
\end{defn}

Having defined these special elements, we give a few examples.

\begin{example}
\label{example:f_divgroups}

	\begin{itemize}
		\item Let $G$ be an infinite hyperbolic group acting on its Cayley graph with respect to a finite generating set. Then there exists an exponential function $f$ such that $G$ contains an $f$-divergent element (see \cite[Lemma 3.4]{GS21}).
		\item Let $(G,\mathcal P)$ be a relatively hyperbolic group, where all $P\in\mathcal P$ have infinite index in $G$. Then there exists an exponential function $f$ such that $G$ contains an $f$-divergent element for the action on the coned-off Cayley graph $\Cay(G,S\cup \mathcal P)$, where $S$ is a finite generating set for $G$ (see [\cite[Lemma 3.5]{GS21}).
  \item Let $G$ be an acylindrically hyperbolic group, acting acylindrically and non-elementarily on the hyperbolic space $X$. Then any loxodromic element $g$ for this action is $f$-divergent for some diverging function $f$ (depending on $g$). This essentially follows from \cite[Proposition 10.4]{MathieuSisto} as we show below in Lemma \ref{lem:fdivelement_acylindrically_hyperbolic} (we believe this also follows from arguments in \cite{SistoHyperboliccalyembedded}).

   \item Let $G$ be the mapping class group of a finite-type surface $S$ (with finitely many exceptions), acting on its curve graph $\mathcal CS$. Then pseudo-Anosov elements are $f$-divergent elements for some linear function $f$; for short, we will say linear-divergent. In fact, a similar result holds for hierarchically hyperbolic groups, see Lemma \ref{lem:HHS_divergent_element} below.
   \item Let $G$ be a group of isometries of a proper CAT(0)-space containing a rank-one element $g_0$. Then $g_0$ is a linear-divergent element, see Lemma \ref{lem:rank_one_divergent} below.
   \end{itemize}
\end{example}

	

\subsection{Proofs of $f$-divergence}

In this subsection we show that the elements described in the last two examples of Example \ref{example:f_divgroups} are indeed $f$-divergent for the claimed $f$. We start with the following well-known result in hyperbolic geometry, which we will use throughout this paper.

\begin{lemma}
\label{lem:exo_hyperbolicspace}
Let $X$ be a $\delta$-hyperbolic space. Let $Q$ be a quasi-convex proper subspace and $\pi_Q:X \rightarrow Q$ a closest-point projection.
Then there exists a constant $R>0$ only depending on $\delta$ and the quasi-convexity constant of $Q$ such that the following holds:
\begin{itemize}
    \item $\pi_Q$ is $R$-coarsely Lipschitz.
    \item For all $x,y \in X$, if $d_{X}(\pi_Q(x), \pi_Q(y)) \geq R$ then there are points $m_1, m_2 \in [x,y]$ such that $d_X(m_1,\pi_Q(x)) \leq R$ and $d_X(m_2,\pi_Q(y)) \leq R$   where $[x,y]$ is a geodesic between $x$ and $y$. Further, the subgeodesic of $[x,y]$ between $m_1$ and $m_2$ lies in the $R$-neighborhood of $Q$.
\end{itemize}
\end{lemma}

As was stated above, the following lemma will follow from \cite[Proposition 10.4]{MathieuSisto}.

\begin{lemma}
\label{lem:fdivelement_acylindrically_hyperbolic}

Let $G$ act acylindrically and non-elementarily on the hyperbolic space $X$. Then any loxodromic element $g$ for this action is $f$-divergent for some diverging function $f$ (depending on $g$).
 \end{lemma}
 
 \begin{proof}
 
 We note that the statement of \cite[Proposition 10.4]{MathieuSisto} involves some 'acylindrically intermediate'  space $Y$ for for $(G,X)$; by \cite[Section 10]{MathieuSisto} one can take $Y=G$ (with a word metric), which is what we will do in this proof.
 
 	Let $g$ be a loxodromic element as in the statement of the lemma, where the quasi-isometric embedding $\langle g \rangle \hookrightarrow X $ has constants $(a,b)$. Then a path along $\langle g\rangle $ is a $L$-Lipschitz discrete path for some $L$ depending only on $a$ and $b$. Let $F$ be the constant from \cite[Proposition 10.4]{MathieuSisto} associated to $L$ and let $\rho : \mathbb R^{+} \to \mathbb R^{+}$ be the corresponding diverging function. Let $d'$ be such that for all $d\geq d'$, we have $\rho(d) \geq 2(a+1).$ Let $\theta \geq a(8R+2F+b)$ where $R$ is from Lemma \ref{lem:exo_hyperbolicspace}, and let: \begin{equation*}
    f(d) =
    \begin{cases*}
      0  & if $d<d'$ \\
     \frac{\rho(d)\theta}{2}        & if $d \geq d'$.
    \end{cases*}
  \end{equation*}
 	Let $\alpha$ be a $L$-Lipschitz (discrete) path, say from $\alpha^{-}$ to $\alpha^{+}$), such that $d_G(\pi(\alpha^{-}), \pi(\alpha^{+})) \geq \theta$ and staying outside of the $d$-neighborhood of $\langle g \rangle$. We want to show that $\ell (\alpha) >f(d)$.
 	
 	Let $\beta$ be the discrete path along $\langle g \rangle$ from  $\pi(\alpha^{-})$ to $\pi(\alpha^{+}$. Then, by \cite[Proposition 10.4]{MathieuSisto}, as $\beta$ is $L$-Lipschitz, we have:
 $$ \max\{ \ell_G(\alpha), \ell_G(\beta)\} \geq \left( d_X(\alpha^{-}, \alpha^{+})-d_X(\alpha^{-}, \pi(\alpha^{-}))-d_X(\alpha^{+}, \pi(\alpha^{+}))-F\right)\rho(d).$$
 Now, by Lemma \ref{lem:exo_hyperbolicspace} there exists $R>0$ only depending on $\langle g \rangle$ such that: $$d_X(\alpha^{-}, \alpha^{+})-d_X(\alpha^{-}, \pi(\alpha^{-}))-d_X(\alpha^{+}, \pi(\alpha^{+}) \geq d_X(\pi(\alpha^{-}), \pi(\alpha^{+}))-4R.$$
 
 If $d<d'$, then $f(d)=0$ and there is nothing to show. Hence, we only need to consider $d>d'$. By the choice of $\theta$, we have that $d_X(\pi(\alpha^{-}), \pi(\alpha^{+})) \geq 4R+F$, hence for all $d>d'$: \begin{align*}
 \begin{split}
 	\ell_G( \beta) &\leq ad_X(\pi(\alpha^{-}), \pi(\alpha^{+}))+b<(a+1)d_X(\pi(\alpha^{-}), \pi(\alpha^{+}))\\
 	& \leq \frac{\rho(d)}{2} d_X(\pi(\alpha^{-}), \pi(\alpha^{+}))<\rho(d)\left(d_X(\pi(\alpha^{-}), \pi(\alpha^{+}))-4R-F\right).
 \end{split}	
 \end{align*}

 Hence, $\ell_G(\alpha) \geq \rho(d)\left(d_X(\pi(\alpha^{-}), \pi(\alpha^{+}))-4R-F\right) > \rho(d)\theta/2=f(d)$, as required.
 	 \end{proof}

 The following lemma is similar to what was done in the proof of \cite[Proposition 5.9]{TrackingratesRW} for the mapping class group, which only relies on the HHS structure of $G$; we include the lemma here for completeness. We refer the reader to \cite{HHS1, HHS2, WhatisHHS} for general background on hierarchically hyperbolic spaces, which we do not recall in detail here since in this paper hierarchical hyperbolicity is only used in the following lemma. In short, we will use that a hierarchically hyperbolic structure $(X, \mathfrak S)$ gives a collection of hyperbolic spaces $\mathcal C Y$, for $Y\in \mathfrak S$, and there is a ``main'' (formally, $\nest$-maximal) hyperbolic space $\mathcal C S$. Moreover, there are coarsely lipschitz (projection) maps $\pi_Y:X\to\mathcal C Y$ and we will use mostly two things. First, the distance formula, which gives an estimate on the distance between two points in $X$ in terms of their projections to the various hyperbolic spaces; we spell out the part of the estimate we need in the proof below. Second, we need the Bounded Geodesic Image axiom, which says that if $x,y$ project far in some $\mathcal C Y$, for $Y\neq S$, then any geodesic from $\pi_S(x)$ to $\pi_S(y)$ passes uniformly close to a specified point $\rho^Y_S$.

\begin{lemma}
\label{lem:HHS_divergent_element}
	Let $G$ be an hierarchically hyperbolic group  (with structure $(G, \mathfrak S)$) and with an element $g_0$ acting loxodromically on the main hyperbolic space $\mathcal CS$. Then $g_0$ is a linear-divergent element for this action.
\end{lemma}

\begin{proof}
We let $\pi: G \to \langle g_0 \rangle$ be a $\mathcal CS$-projection, as defined in Definition \ref{defn:X_proj}, where $S$ is the $\nest$-maximal element of $\mathfrak S$. Let $L_0$ such that $L_0/2 \geq \max\{12E, 2R, s_0\}$, where $E$ is the constant from the Bounded Geodesic Image axiom, $R$ is from Lemma \ref{lem:exo_hyperbolicspace} for the hyperbolic space $\mathcal CS$ and the quasi-convex subset $\langle g_0\rangle x_0$ and $s_0$ is from the distance formula. Define $f(d):=(\lambda_2^{-1}\lambda_1^{-1}(d-\nu_1)-\nu_2)/2$ where $(\lambda_1, \nu_1)$ are the constants associated to $2L_0$  and $(\lambda_2, \nu_2)$ are the constants associated to $L_0$ in the distance formula. Let $\theta$ be such that $\theta/a-b \geq 20E $ where $(a,b)$ are the quasi-geodesic constant for $\langle g_0 \rangle x_0 \subseteq \mathcal CS$. We note that both $f$ and $\theta$ only depend on the HHS structure of $G$ and on $g_0$.
 
  Let $\alpha$ be a path in $G$ as in Definition \ref{defn:H_T}, from $\alpha^{-} $ to $\alpha^{+}$, avoiding the $d$-neighborhood of $\langle g_0 \rangle$ and such that $D:=d_G(\pi(\alpha^{-}), \pi(\alpha^{+})) \geq \theta$.

	The main ingredient to this proof is to show that if a domain $U \in \mathfrak S$ is such that $d_U(\alpha^{-}, \pi_{\gamma}(\alpha^{-}))$ is big then so is $d_U( \alpha^{-}, \alpha^{+})$. For two points $x,y \in G$ and $L \geq 0$, we say that $U \in \textit{REL}(x,y,L)$ if $d_U(x,y) \geq L$.
	
	\begin{claim}
	\label{claim:relevant}
For all $L \geq L_0$, if $d_U(\alpha^{-}, \pi(\alpha^{-})) \geq L$ then $d_U(\alpha^{-}, \alpha^{+}) \geq d_U(\alpha^{-}, \pi(\alpha^{-}))/2$. In particular, if $U \in \textit{REL}(\alpha^{-}, \pi(\alpha^{-}), L)$, then $U \in \textit{REL}(\alpha^{-}, \alpha^{+}, L/2).$
	\end{claim}
	
	\begin{proof}[Proof of Claim]
	We first note that there exists a constant $K_1$ such that for all $U \in \mathfrak S$ the projection map $\pi_U:G \to \mathcal CU$ is $(K_1, K_1)$-coarsely Lipschitz. Hence $d_U(\pi(\alpha^{-}), \pi(\alpha^{+})) \leq K_1D+K_1.$  let $L \geq L_0$ and let $U \in \textit{REL}(\alpha^{-}, \pi(\alpha^{-}), L)$.
	
	The first case is if $U \neq S$, then by the Bounded Geodesic Image axiom,  there exists a constant $E$ such that $\rho^{U}_S \cap \mathcal N_{E}([\pi_S(\alpha^{-}), \pi_S(\pi(\alpha^{-})]) \neq \emptyset$.
	By hyperbolicity of $\mathcal CS$ and the choice of $\theta$, we have that $\rho^{U}_S \cap \mathcal N_{E}([\pi_S(\alpha^{+}), \pi_S(\pi(\alpha^{+})]) = \emptyset$. Hence, using the Bounded Geodesic Image axiom, we can bound $d_U(\alpha^{+}, \pi(\alpha^{+})) \leq E$.
	Now, let $\beta \subseteq [\pi_S(\pi(\alpha^{-})), \pi_S(\pi(\alpha^{+}))]$ be the points $x$ such that $d_S(x, \pi(\alpha^{-})) \geq 10E$. Then we have that $\rho^{U}_S \cap \mathcal N_E(\beta) =\emptyset$ and so by the Bounded Geodesic Image axiom, we get that $d_U(\beta^{-}, \beta^{+}) \leq E$. Noting that $d_U(\pi(\alpha^{-}), \pi(\alpha^{+}))=d_U(\beta^{-}, \beta^{+})+10E \leq 11E$, we get that:
	$$ d_U(\alpha^{-}, \alpha^{+}) \geq d_U(\alpha^{-}, \pi(\alpha^{-}))-d_U(\pi(\alpha^{-}), \pi(\alpha^{+}))-d_U(\pi(\alpha^{+}), \alpha^{+}) \geq d_U(\alpha^{-}, \pi(\alpha^{-}))-12E\geq d_U(\alpha^{-}, \pi(\alpha^{-}))/2. $$ 
	
	If $U =S$, then by hyperbolicity of $\mathcal CS$ and Lemma \ref{lem:exo_hyperbolicspace} we have that: $$ d_S(\alpha^{-}, \alpha^{+}) \geq d_S(\alpha^{-}, \pi(\alpha^{-}))-2R \geq d_S(\alpha^{-}, \pi(\alpha^{-}))/2. $$
	
	This finishes the proof of the claim.
	\end{proof}
	Let $L=2L_0$, then by the choice of $(\lambda_1,\nu_1)$ above the claim, and the distance formula we have that: $$ d \leq d_G(\alpha^{-}, \pi(\alpha^{-})) \leq \lambda_1\sum_{U \in \mathfrak S} \ignore{d_U(\alpha^{-}, \pi(\alpha^{-}))}{L} +\nu_1.   $$ 
	
	Where we recall that $\ignore{A}{L}=A$ if $A \geq L$ and $\ignore{A}{L}=0$ otherwise.
	
	By the choice of constants $(\lambda_2, \nu_2)$ for $L_0$ we have: 
	\begin{align*}
		\begin{split}
		 d_G(\alpha^{-}, \alpha^{+}) &\geq \lambda_2^{-1}\sum_{U \in \mathfrak S} \ignore{d_U(\alpha^{-}, \alpha^{+})}{L/2} -\nu_2 \\
			&\geq \lambda_2^{-1}\sum_{U \in \textit{REL}(\alpha^{-}, \pi(\alpha^{-}), L)} \ignore{ d_U(\alpha^{-}, \pi(\alpha^{-}))}{L}-\nu_2 \\
			&\geq \lambda_2^{-1}\lambda_1^{-1}(d-\nu_1)-\nu_2 >f(d)
		\end{split}
	\end{align*}
	where we go from the first to the second line by using Claim \ref{claim:relevant}. Hence $\ell(\alpha) \geq d_G(\alpha^{-}, \alpha^{+}) >f(d). $
\end{proof}

The following lemma shows the existence of a linear-divergent element for CAT(0) groups with a rank-one isometry. We do this by looking at the associated action of the group on the BBF quasi-tree and using the contracting properties of this element. It might be possible to use a similar argument when in place of using the BBF quasi-tree one uses the hyperbolic space defined for a CAT(0) space in  
\cite{SprianoPetytZalloum}. 
\begin{lemma}
\label{lem:rank_one_divergent}
	Let $G$ be a group of isometries of a proper CAT(0) space with a rank-one element $g_0$. Then $g_0$ is a linear-divergent element.
\end{lemma}
\begin{proof}
	We note that by \cite[Examples 2.1-(3), Theorem H]{BBF}, there is a quasi-tree $X$ such that $G$ has a non-elementary action on $X$, and $g_0$ is loxodromic for this action.
	
	Further, the $X$-projection $\pi:G \to \langle g_0 \rangle$ is in fact within bounded error of a closest point projection, see  \cite[Lemma 4.9]{BBF} (and the comment above said lemma). Now, $g_0$ is loxodromic for the action on $\mathbb Y$ and hence by \cite[Theorem 5.4]{BestvinaFujiwaraCohomology}, is $B_1$-contracting for the closest point projection $\pi :G \to \langle g_0 \rangle$. We show that this implies that $g_0$ is a linear divergent element. 
	
	Indeed, let $\alpha$ be a path from $\alpha^{-}$ to $\alpha^{+}$ staying outside of the $d$-neighbourhood of $\langle g_0 \rangle$ and such that the $\mathbb Y$-projection (and hence the closest point projection $\pi:G \to \langle g_0 \rangle$) satisfies $d\left(\pi(\alpha^{-}), \pi(\alpha^{+})\right) \geq \theta:=10B_1$ where $B_1$ is from the axis of $g_0$ being $B_1$-contracting.
	
	We first note that $B(\alpha^{-}, d/2) \cap \langle g_0 \rangle = \emptyset$, hence $\diam \left(\pi(B(\alpha^{-}, d/2)\right)) \leq B_1$ and so in particular $\alpha $ is not properly contained in $B(\alpha^{-}, d/2)$ and hence $\ell (\alpha) \geq d/2$, as required. 
\end{proof}

\subsection{Divergence and WPD elements}

In this subsection we recall the notion of WPD element and we argue that $f$-divergent elements, for $f$ a diverging function, are WPD. As we recall WPD elements, we also establish some notation, following \cite{GS21}, that will be used in the statements of our main results and throughout this paper. Finally, we establish a connection between $f$-divergence and divergence of groups.

\subsubsection{WPD elements}

We recall the definition of a WPD element from \cite{bestvinafujiwara}.

\begin{defn}[WPD element]
Let $G$ be a group acting on a hyperbolic space $X$ and $g$ an element of $G$. We say that $g$ satisfies the \textit{weak proper discontinuity condition} (or that $g$ is a WPD element) if for all $\kappa >0$ and $x_0 \in X$ there exists $N \in \mathbb{N}$ such that \[
\# \{ h \in G \vert \quad d_X(x_0, hx_0)< \kappa, \quad d_X(g^{N}x_0, hg^{N}x_0) < \kappa \} < \infty.
\]
\end{defn}

The following lemma shows that whenever $G$ contains an $f$-divergent element $g_0$ then $g_0$ is a loxodromic WPD element.  

\begin{lemma}
\label{lem:divergent_element_is_WPD}
If $g_0$ is an $f$-divergent element (for an increasing divergent function $f$) for some action on a hyperbolic space $X$ then $g_0$ is a loxodromic WPD element for this action on $X$.
\end{lemma}

\begin{proof}
	Identical to the proof of \cite[Lemma 3.8]{GS21}. In the set-up of that paper, $f$ is assumed to be super-linear, but for the proof of that lemma it only needs to be divergent.
\end{proof}

Recall that any loxodromic WPD element $g$ is contained in a unique maximal elementary subgroup of $G$, denoted $E(g)$ and called the \textit{elementary closure} of $g$, see  \cite[Theorem 1.4]{Osin_acylindrical}.

We will refer to the following lemma as the strong Behrstock inequality, and we will use it very often. The lemma follows from \cite[Theorem 4.1]{bromberg2017acylindrical}, see \cite[Section 3.3]{GS21}.

We briefly explain the meaning of the lemma. Given a point $x$ and two distinct cosets of $E(g)$ for a loxodromic element $g$, (one version of) the Behrstock inequality says that on one of the two cosets the $X$-projection of $x$ is close to the $X$-projection of the other coset. The strong version says that, up to perturbing the $X$-projections, we can actually ensure that, rather than being close, the projections actually coincide.

 \begin{lemma}
\label{lem:Behrstock}
Let $g$ be a loxodromic WPD element. Then, for $\gamma=E(g)$, there is a $g$-equivariant map $\pi_{\gamma}:G\to \mathcal P(\gamma)$, where $\mathcal P(\gamma)$ is the set of all subsets of $\gamma$, and a constant $B$ with the following property. For all $x\in G$ and $h\gamma \neq h'\gamma$ we have
\[ d_X(\pi_{h\gamma}(x),\pi_{h\gamma}(h'\gamma)) >B \implies \pi_{h'\gamma}(x)=\pi_{h'\gamma}(h\gamma),
\]
where $\pi_{k\gamma}(z)=k\pi_{\gamma}(k^{-1}z)$. Moreover, for all $x\in G$ the Hausdorff distance between $\pi_\gamma(x)$ and an $X$-projection of $x$ to $\langle g\rangle$ is bounded by $B$.
\end{lemma}

As in \cite[Section 3.3]{GS21}, we will often consider projections distances:

\begin{notation}
\label{not:d_hgamma}
For a loxodromic WPD element $g$, consider $\gamma=E(g)$ and the maps $\pi_{h\gamma}$ as in Lemma \ref{lem:Behrstock}. For $h,x,y\in G$ we denote
$$d_{h\gamma}(x,y)=\diam (\pi_{h\gamma}(x)\cup\pi_{h\gamma}(y)),$$
and similarly when we project subsets rather than points.
\end{notation}



As in the statement of Lemma \ref{lem:Behrstock}, we will often denote $\gamma=E(g_0)$ when a loxodromic WPD element $g_0$ has been fixed. We  often look at the set of cosets where two given elements have far away projections, as captured by the following definition.

 \begin{defn}
 \label{defn:H_T}
 Given a loxodromic WPD element $g$, for $x,y\in G$ and $T\geq 0$ we define the set \[ \mathcal{H}_T(x,y) :=\{ h\gamma \hspace{1mm} : \hspace{1mm} d_{h\gamma}(x,y) \geq T\},
\]
and similarly when we projects sets rather than points.
 \end{defn}
 
 \subsubsection{Divergence}
 
There are many definitions of divergence for a metric space. These different definitions are equivalent under mild conditions, see \cite{behrstockDrutu}, \cite{DrutuMozesSapir} for more on this. More precisely, the various definition yield functions that are equivalent with respect to the equivalence relation generated by the partial order where $f\preceq g$ if there exists a constant $C$ such that $f(x)\leq Cg(Cx+C)+C$. We will use the definition given in \cite[Definition 3.1]{behrstockDrutu}.
 
 \begin{defn}
 \label{defn:divergence}
 	 Consider a constant $0<\delta<1$. For a triple of points $a,b,c \in X$ with $dist(c, \{a,b\}) =r>0$ we define $\dive(a,b,c,\delta)$ to be the infimum of the lengths of all paths connecting $a$ to $b$ and avoiding the ball $B(c, \delta r)$. If no such path exists, we define $\dive(a,b,c,\delta)=\infty$.
 \end{defn}

 The divergence function (given a fixed parameter) of a metric space is then defined considering the supremum over all pairs of points $a,b$ within a given distance of each other; we will not need the precise definition.

As noted in the introduction, the notion of $f$-divergence is used implicitly in the literature to give lower bounds on divergence. We will not need the following lemma, but we point it out to make the connection between $f$-divergence and divergence explicit. Lemma \ref{lem:argument_divergence_path_ouside_ball} below is similar in spirit.

\begin{lemma}
\label{lem:lower_bound_divergence}
	If the group $G$ contains an $f$-divergent element, then  it has divergence $\succeq n f(n)$.
\end{lemma}

\begin{proof}

The proof is identical to \cite[Lemma 3.6]{GS21}, where $f$ was assumed to be superlinear and the conclusion is that the divergence is superquadratic. However, the last displayed formula in the proof gives the bound we require here.
\end{proof}

\section{Statements of main results}

\label{sec:statements}

\subsection{Assumptions on Markov chains}
Throughout this paper, we will assume that our Markov chains have the following properties. We note that the non-amenability condition from \cite[Definition 2.5]{GS21} is not required here.

We use the notation of \cite[Section 2]{GS21}, in particular $w_n^o$ denotes an $n$-step Markov chain starting at $o$, and transition probabilities are denoted $p(\cdot,\cdot)$.

\begin{Assumpt}
\label{assumpt:weak_tameness}
Let the Markov chain $(w_n^{o})$ on $G$ satisfy the following conditions.
	\begin{enumerate}
    \item\label{item:bounded_jumps} {\bf (Bounded jumps)} There exists a finite set $S\subseteq G$ such that $p(g,h)=\mathbb P[w^g_1=h]=0$ if $h\notin gS$.

  \item\label{item:irred} {\bf (Irreducibility)} For all $s \in G$ there exist constants $\epsilon_s, K_s>0$ such that for all $g \in G$ we have
  $$\mathbb P[w^g_k=gs] \geq \epsilon_s$$
  for some $k \leq K_s$.
	\end{enumerate}
\end{Assumpt}

 \begin{remark}
 \label{rmk:bounded_jumps}
 	Once we have fixed a word metric $d_G$ on $G$ then the assumption of Bounded jumps \ref{assumpt:weak_tameness}-\eqref{item:bounded_jumps} is equivalent to the following: There exists a constant $K>0$ such that for all $n \in \mathbb N$ and starting point $o \in G$ we have $d_G(w^o_n, w^o_{n+1}) \leq K$.
\end{remark}

We now state some assumptions on the groups we will consider. We note that these assumptions on $G$ will be different depending on whether we are considering a random walk or a more genral Markov chain. For simplicity of notation, we write $(w_n^{o})$ for a random walk or a Markov chain on $G$. In the case of a random walk, the starting point (which will be considered to be the identity) is not as relevant as in the case of a Markov chain.  In the case of the Markov chain, we will further require that $(w_n^{o})_n$ satisfies the conclusion of \cite[Proposition 5.1]{GS21}. We note that this proposition is enough to show linear progress, in $X$, of the Markov chain, see \cite[Section 6]{GS21}. Here are the assumptions under which we will work:

\begin{assumptx}
\label{assump:random_walk}
$G$ is a finitely generated group containing an $f$-divergent element $g_0$ and $(w^o_n)_n$ is a random walk driven by a finitely supported measure that generates $G$ as a semigroup.
\end{assumptx}

\begin{assumptx}
\label{assump:markov_chain}
$G$ is a finitely generated group containing an $f$-divergent element $g_0$, and $(w_n^{o})_n$ is a Markov chain on $G$ satisfying Assumption \ref{assumpt:weak_tameness} and satisfying the conclusion of \cite[Proposition 5.1]{GS21}, which we recall later.
\end{assumptx}

Various groups are shown in \cite{GS21} to satisfy Assumption \ref{assump:markov_chain} for \emph{tame} Markov chains (\cite[Definition 2.5]{GS21}). These include:

\begin{itemize}
	\item Non-elementary hyperbolic and relatively hyperbolic groups. 
	\item Acylindrically hyperbolic $3$-manifold groups.
\end{itemize}

In work in progress of Mark Hagen, Harry Petyt, and the authors, we show that suitable Markov chains on hierarchically hyperbolic groups also satisfy Assumption \ref{assump:markov_chain}.

\subsection{Main results}

We fix once and for all the $f$-divergent element $g_0$ which is loxodromic and WPD for an action of $G$ on a hyperbolic space $X$. We also fix $\gamma=E(g_0)$ to be the elementary closure of $g_0$. 
The following results will later be used as assumptions for Theorem \ref{thm:lower_bound} below on ``random divergence'' and can be used as a ``black box'', meaning that if one can establish the conclusion of these for other groups, then the theorem will also hold in this case.

The following proposition states that the number of cosets where the distance between the projection of $p$ and $w_n^p$ is big grows linearly. This result is perhaps what the authors should have established in \cite{GS21} in place of their Proposition 5.1. The difference between these propositions is that we want the number of cosets to grow linearly here, whilst in \cite{GS21}, we wanted the sum of the projection distances over the cosets to grow linearly.

\begin{prop}
\label{prop:Axiom_blackbox}	
Let $G$ and $(w^{o}_n)_n$ satisfy either Assumption \ref{assump:random_walk} or Assumption \ref{assump:markov_chain}. Then there exists a constant $T_0$ such that the following holds. For all $T \geq T_0$, there exist  constants $C_0, \epsilon_0>0$ such that for all $p \in G$ and $n \in \mathbb N$ we have:
$$\mathbb P \Big[ \vert \mathcal H_T(p,w_{n}^p) \vert \geq \epsilon_0 n \Big] \geq 1-C_0e^{-n/C_0}. $$
\end{prop}

The second result will be an easy consequence of the conclusion of \cite[Proposition 5.1]{GS21} and roughly states that given a set $\mathcal H_T(o,p)$ of cosets, the probability that a Markov chain starting at $p$ undoes the projections on more than $t$ of these cosets decays exponentially in $t$.

\begin{prop}
	\label{prop:axiom_implies_small_interesction}
	Let $G$ and $(w^{o}_n)_n$ satisfy either Assumption \ref{assump:random_walk} or Assumption \ref{assump:markov_chain}.  Then there exists $T_0'\geq 0$ with the following property. For all $T \geq T_0'$ there exists a constant $C_1>0$ such that the following holds for all $o,p \in G$ and for all $m$:
	
	$$\mathbb P \Big[\vert  \mathcal H_T(o,p) \cap \mathcal H_T(p,w_m^p) \vert > t  \Big] \leq C_1e^{-t/C_1}. $$
\end{prop}

Finally, theorem below states that, under our assumptions, the lower bound on divergence given in Lemma \ref{lem:lower_bound_divergence} is true generically:

\begin{thm}
\label{thm:lower_bound}
	Let $G$ be a group acting non-elementarily on a hyperbolic space $X$, and let $g_0$ be an $f$-divergent element for the action. Let $(w^{p}_n)_n, (z_n^p)_n$ be two independent copies of a Markov chain (or random walk) satisfying the conclusions of both Proposition \ref{prop:Axiom_blackbox} and Proposition \ref{prop:axiom_implies_small_interesction}. Then there exists a constant $\delta_0 \in (0, 1)$ such that for all $\delta \leq \delta_0$, there  exists a constant $C>0$ such that:

$$\mathbb P\Big[\dive(w^{p}_n,z^{p}_n,p, \delta) > n f(n/C)/C\Big]\geq 1-C e^{-n/C}$$
\end{thm}


The rest of the paper is devoted to the proofs of the propositions and the theorem.

\begin{notation}
    From now on, $G$ will always denote a group acting non-elementarily on a hyperbolic space $X$, and $g_0$ will be an $f$-divergent element for the action. We will use all the notations set in Section \ref{sec:f_div}, such as $E(g_0)$, etc.
\end{notation}

\subsection{Some more set up}

Before moving on to the proofs of the propositions and the theorem, we recall some more results about hyperbolic spaces and WPD elements which we will use throughout.  

 Let $x,y \in G$. In view of the fact that the projections on the $h\gamma$ satisfy the projection axioms of \cite{BBF}, by \cite[Theorem 3.3 (G)]{BBF}, we have a total order on $\mathcal{H}_T(x,y)$, as specified below. Regarding the equivalence of the various conditions, this follows from Lemma \ref{lem:Behrstock} (Behrstock inequality).

\begin{lemma}(Consequence of \cite[Theorem 3.3 (G)]{BBF})
\label{lem:linearorder}
Fix a loxodromic WPD element $g$. For any sufficiently large $T$ and for any $x,y \in G$, the set $\mathcal{H}_T(x,y) \cup \{x,y\}$ is totally ordered with least element $x$ and greatest element $y$. The order is given by $h\gamma \prec h'\gamma$ if one (and hence all) of the following equivalent conditions hold for $B$ as in Lemma \ref{lem:Behrstock}:
\begin{itemize}
    \item $d_{h\gamma}(x,h'\gamma) > B$.
    \item $\pi_{h'\gamma}(x)= \pi_{h'\gamma}(h\gamma)$.
    \item $d_{h'\gamma}(y,h\gamma) > B$.
    \item $\pi_{h\gamma}(y)= \pi_{h\gamma}(h'\gamma)$.
\end{itemize}
\end{lemma}

We now study how various sets of large projections can differ, starting with the following lemma.

 \begin{lemma}(\cite[Lemma 3.15]{GS21})
\label{lem:2cosetsmax}
Fix a loxodromic WPD element $g$. For any $T\geq 10 B$ which satisfies Lemma \ref{lem:linearorder}, with $B$ as in Lemma \ref{lem:Behrstock}, the following holds.
Let $x,y,z\in G$. Then there at most $2$ cosets $h\gamma \in \mathcal{H}_T(x,y)$ such that $\pi_{h\gamma}(z)$ is distinct from both $\pi_{h\gamma}(x)$ and $\pi_{h\gamma}(y)$.
\end{lemma}

A consequence of the lemma is given below. We shall use this lemma several times and in the two slightly different forms given in the statement. We include proofs for completeness, but we note that they could be deduced from \cite[Lemma 3.6]{bromberg2017acylindrical}.

\begin{lemma}
\label{lem:bound_number_cosets}
	For all $x,y,z \in G$, there are at most $2$ cosets $t_1\gamma, t_2\gamma \in \mathcal H_T(x,z)$ such that $$ \mathcal H_T(x,z) \backslash \mathcal H_T(z,y) \subseteq \mathcal H_T(x,y) \cup \{t_1\gamma, t_2\gamma \}.  $$
In particular:
	$$ \vert \mathcal H_T(x,y) \vert \geq \vert \mathcal H_T(x,z) \backslash \mathcal H_T(z,y)\vert -2\geq \vert  \mathcal H_T(x,y)\backslash \mathcal H_T(z,y) \vert-4. $$
\end{lemma}

\begin{proof}
	By Lemma \ref{lem:2cosetsmax}, there are at most 2 cosets (call them $t_1 \gamma $ and $t_2\gamma)$ in $\mathcal H_T(x,z)\backslash \mathcal H_T(z,y)$ such that $\pi_{t_i\gamma}(y) \neq \pi_{t_i\gamma}(x)$ and  $\pi_{t_i\gamma}(y) \neq \pi_{t_i\gamma}(z)$. Now for all  $h \gamma \in \mathcal H_T(x,z)\backslash \mathcal H_T(z,y)$ but these two $t_i \gamma$  we must have $\pi_{h\gamma}(y)=\pi_{h\gamma}(z)$ (otherwise $h\gamma \in  \mathcal H_T(z,y)$). Hence for all but the two $t_i \gamma$ cosets, we have $d_{h\gamma}(x,y) =d_{h\gamma}(x,z) \geq T$ as $h\gamma \in \mathcal H_T(x,z)$ and so all but at most 2 cosets of $\mathcal H_T(x,z)\backslash \mathcal H_T(z,y)$ belong to $\mathcal H_T(x,y)$. This proves the first part of the lemma and thus the first inequality. Swapping $z$ and $y$, we see that $\mathcal H_T(x,y)\backslash \mathcal H_T(z,y)$ is contained in $\mathcal H_T(x,z)$ possibly up to two cosets. But then the same is true with $\mathcal H_T(x,z)\backslash \mathcal H_T(z,y)$ replacing $\mathcal H_T(x,z)$, showing the second inequality.
\end{proof}

The following states that we have a 'coarse triangular inequality' for the number of cosets in $\mathcal H_T$.

\begin{cor}
\label{cor:triangle_inequality}
	For all $x,y,z \in G$ we have 
	
	$$ \vert \mathcal H_T(x,z) \vert \leq  \vert \mathcal H_T(x,y) \vert+\vert \mathcal H_T(y,z) \vert+2.$$
	
%
%
\end{cor}
\begin{proof}
	This follows immediately from Lemma \ref{lem:bound_number_cosets} and the fact that $$\vert \mathcal H_T(x,z) \backslash \mathcal H_T(z,y)\vert =\vert \mathcal H_T(x,z)\vert - \vert  \mathcal H_T(x,z) \cap \mathcal H_T(z,y) \vert \geq \vert \mathcal H_T(x,z)\vert - \vert  \mathcal H_T(z,y) \vert. $$ 
\end{proof}

\begin{notation}
	Once a finite generating set has been fixed for $G$, we fix a word metric $d_G$ on $G$ and note that for an action  $G \acts X$ for all $g,h\in G$ we have $d_X(gx_0,hx_0) \leq d_G(g,h)$ up to rescaling $X$, which we do from now on for convenience.
\end{notation}

The following lemma says roughly that $X$-projections are closest-point projections up to multiplicative and additive error.

 \begin{lemma}
\label{lem:distances_projection}
	There exists a constant $D$ (depending only on $\gamma$) such that for all $h\gamma$ and $p \in G$, we have $$ d_G(p, \pi_{h\gamma}(p)) \leq Dd_G(p,h\gamma)+D$$
\end{lemma}

\begin{proof}
We let $p' \in h\gamma$ be such that $d_G(p,p')=d_G(p, h\gamma)$. Since $\pi_{h\gamma}$ is at bounded distance from a closest-point projection (Lemma \ref{lem:Behrstock}), the conclusion of Lemma \ref{lem:exo_hyperbolicspace} holds for $\pi_{h\gamma}$; we denote by $R$ the corresponding constant. Also, let $a,b$ be such that $h\gamma$ is quasi-isometrically embedded in $X$ with constants $(a,b)$ (meaning that $h\gamma$ endowed with the restriction of $d_G$ is $(a,b)$-quasi-isometrically embedded in $X$ via the orbit map for $x_0$). We set $D=\max\{1+a, a(b+R+2B)\}$.
If $d_X(p'x_0, \pi_{h\gamma}(p)x_0) \leq R$ we have that $d_G(p',\pi_{h\gamma}(p)) \leq aR+ab$ and so $d_G(p, \pi_{h\gamma}(p)) \leq d_G(p,p')+aR+ab \leq Dd_G(p,h\gamma)+D.$
If $d_X(p'x_0, \pi_{h\gamma}(p)x_0) > R$, then by Lemma \ref{lem:exo_hyperbolicspace} there exists a point $t \in [px_0,p'x_0]$ such that $d_X(t, \pi_{h\gamma}(p)x_0) \leq R$ and so
$$d_G(p,p') \geq d_X(px_0, p'x_0) \geq  d_X(px_0, t) +d_X(p'x_0, \pi_{h\gamma}(p)x_0)-R-2B \geq \frac{1}{a}d_G(p', \pi_{h\gamma}(p))-b-R-2B,$$
where the "$2B$" is due to the diameter of $\pi_{h\gamma}(p)$. This leads to
$$d_G(p, \pi_{h\gamma}(p)) \leq d_G(p,p')+d_G(p',\pi_{h\gamma}(p))\leq (1+a)d_G(p,p')+a(b+R+2B) \leq Dd_G(p,h\gamma)+D.$$
\end{proof}

\subsection{Fixing constants}

 For a basepoint $x_0 \in X$, as $h\gamma x_0$ is a quasi-convex subset of the hyperbolic space $X$, we have by Lemma \ref{lem:exo_hyperbolicspace} that the closest point projection to $h\gamma x_0$ is coarsely Lipschitz. Now, $h\gamma$ is quasi-isometrically embedded in $X$, hence the $X$-projections $\pi_{h\gamma}$ are $L$-coarsely Lipschitz, for some constant $L \geq 0$ which we fix from now on (see below). We also fix the following constants. Recall that we have already fixed an $f$-divergent element $g_0$ and its elementary closure $\gamma =E(g_0)$.
\begin{notation}
\label{not:all_constants_fixed}
\begin{itemize}
	\item Let $L$ be such that the projections $\pi_{h\gamma}$ are $L$-coarsely Lipschitz and let $R$ be from Lemma \ref{lem:exo_hyperbolicspace} corresponding to the quasi-convex subset $h\gamma x_0 \subseteq X$.
	\item Let $\theta$ and $f$ be from the definition of $g_0$ being an $f$-divergent element.
	\item $B$ is from the strong Behrstock inequality (Lemma \ref{lem:Behrstock}).
	\item Let $K$ be from the bounded jumps condition in Remark \ref{rmk:bounded_jumps}. We note that a sample path of our Markov chain will be a discrete path with $K$-bounded jumps, that is, sequences of points with consecutive ones being distance at most $K$ apart.

\end{itemize}
\end{notation}

\section{Probabilistic arguments}
\label{sec:probability}

We will show that groups satisfying either of Assumptions \ref{assump:random_walk} or \ref{assump:markov_chain} satisfy the conclusion of Proposition \ref{prop:Axiom_blackbox} and Proposition \ref{prop:axiom_implies_small_interesction}. This will then allow us later to prove Theorem \ref{thm:lower_bound}. It might be possible that groups and random walks as in Assumption \ref{assump:random_walk} satisfy the conclusion of \cite[Proposition 5.1]{GS21}, which would give another way of dealing with those.

\subsection{Proof of Proposition \ref{prop:Axiom_blackbox} and Proposition \ref{prop:axiom_implies_small_interesction} for groups satisfying Assumption \ref{assump:random_walk}}

For this subsection, we fix the setup of Assumption \ref{assump:random_walk}.

 We recall some properties of the projection complexes introduced in \cite{BBF}, which we will then use to prove Proposition \ref{prop:Axiom_blackbox}.
 
\subsubsection{Projection complex}
 
We will not define the projection complex formally but rather we recall the properties that are relevant here. Roughly speaking, the projection complex for a loxodromic WPD $g_0$ is made of the cosets of $E(g_0)$, thought of as quasi-lines, with two of them stuck together if there is no third coset on which the two cosets project far away. In particular, the ambient group $G$ is a subset of this projection complex, which we use below; we recall also that $G$ acts on this projection complex. We invite the reader to look at \cite{BBF} or \cite{bromberg2017acylindrical} for details. 

Let $\mathbb Y:=\{ h\gamma :h \in G \}$. By \cite[Example 4.38]{dahmaniguirardelosin}, $\mathbb Y$ satisfies the axioms of \cite[Section 3.1]{BBF}. For every $U>0$, we can then build a projection complex $\mathbb Y_U$ as in \cite[Definition 3.6]{BBF}.

First of all:

\begin{lemma}(\cite[Theorem 3.16]{BBF}) 
\label{lem:existenceU_BBFtree}
For all sufficiently large $U>0$, $\mathbb Y_U$ is a quasi-tree.
\end{lemma}


The reason that projections complexes are useful for us is that the distance between two points is related to the number of cosets on which there is a large proejction: 

\begin{prop}(\cite[ Proposition 3.7]{BBF})
\label{prop:BBF_lowerbound}
	For all sufficiently large $U>0$ and for $x,y$ in $G$, we have $d_{\mathbb Y_U}(x\gamma, y\gamma) \leq \vert \mathcal H_{U}(x,y) \vert +1$.
\end{prop}


We note that we have an action $G \curvearrowright \mathbb Y_U$ on the quasi-tree $\mathbb Y_U$. The following ensures that this action is acylindrical.

\begin{thm}(Special case of \cite[Theorem 1.1]{bromberg2017acylindrical})
\label{thm:action_on_BBF_acylindrical}
	For all sufficiently large $U>0$ the following holds. If $G$ acts on $\mathbb Y_U$ and there exists $N$ such that the common stabiliser of any pair of distinct elements of $\mathbb Y_U$ has cardinality at most $N$, then the action is acylindrical.
\end{thm}

The common stabiliser of two cosets of $E(g_0)$ is the intersection of the corresponding conjugates, which has uniformly bounded size, so that the action of $G$ on $\mathbb Y_U$ is indeed acylindrical.

We will also need the following lemma on the quasi-trees $\mathbb Y_U$, which allows us to control Gromov products in terms of projections.

\begin{lemma}
\label{lem:gromov_proj}
For all sufficiently large $U>0$ there exists $\Delta>0$ such that the following holds. Fix distinct $h_i\gamma$ for $i=1,2,3$ and let $\mathcal H=\mathcal H_U(h_1\gamma,h_3\gamma)\cap \mathcal H_U(h_2\gamma,h_3\gamma)$. Then
$$(h_1\gamma, h_2\gamma)_{h_3\gamma} > |\mathcal H|/2-\Delta.$$
\end{lemma}

\begin{proof}
Throughout the proof we assume that $U$ is large enough that we can apply the results above, as well as other results from \cite{bromberg2017acylindrical}. Since $\mathbb Y_U$ is hyperbolic, the Gromov product in the statement coincides up to an error only depending on the hyperbolicity constant with the distance from $h_3\gamma$ to a geodesic from $h_1\gamma$ to $h_2\gamma$. By \cite[Corollary 3.8]{bromberg2017acylindrical}, said geodesic is at bounded Hausdorff distance from $\mathcal H_U(h_1\gamma,h_2\gamma)$. Hence, we have to give a lower bound on the distance in $\mathbb Y_U$ between $h_3\gamma$ and any element $h\gamma$ of $\mathcal H_U(h_1\gamma,h_2\gamma)$ in terms of $|\mathcal H|$.

By \cite[Lemma 3.6]{bromberg2017acylindrical}, up to replacing $h\gamma$ with an adjacent vertex of $\mathbb Y_U$, and up to swapping $h_1\gamma$ and $h_2\gamma$ we can assume that $h\gamma\in \mathcal H_U(h_3\gamma,h_1\gamma)$. Moreover, that set is totally ordered as in Lemma~\ref{lem:linearorder} again as a consequence of \cite[Theorem 3.3 (G)]{bromberg2017acylindrical}. By \cite[Lemma 3.6]{bromberg2017acylindrical} we have that $h\gamma$ comes after $\mathcal H$ in the order. Therefore, for all $h'\gamma\in\mathcal H$, we have $d_{h'\gamma}(h\gamma,h_3\gamma)\geq U$, that is, $\mathcal H\subseteq \mathcal H_U(h_3\gamma,h\gamma)$. We can then conclude that $d_{\mathbb Y_U}(h\gamma,h_3\gamma)\geq |\mathcal H|/2$ in view of \cite[Corollary 3.7]{bromberg2017acylindrical}, and this concludes the proof.
\end{proof}

\begin{notation}
\label{not:fixing_threshold_T}
Recall that $\theta, B$ and $L$ have been fixed in Notation \ref{not:all_constants_fixed}. Fix $T_0 \geq \max\{100(\theta+B+L), U \}$ large enough that Lemma \ref{lem:linearorder}  (linear order) is satisfied and where $U$ is large enough that all results on $\mathbb Y_U$ stated above apply.
\end{notation}

As in Assumption \ref{assump:random_walk}, we consider a random walk $(w_n)$ on $G$ driven by a finitely supported measure that generates $G$ as a semigroup.

The following is the main theorem of \cite{MaherTiozzo} (see also \cite[Theorem 9.1]{MathieuSisto}):


%

\begin{thm}(\cite{MaherTiozzo})
\label{thm:linear_progress_mt}
Let $G$ be a finitely generated group acting acylindrically on a geodesic metric space $X$. Let $(w_n)_n$ be a random walk driven by a non-elementary finitely supported measure $\mu$. Then there exists a constant $C>0$ such that for all $n$ we have $$ \mathbb P \Big[ d_X(x_0, w_nx_0) \geq n/C \Big] \geq 1-Ce^{-n/C}.$$ 
\end{thm}

By Theorem \ref{thm:action_on_BBF_acylindrical} the action of $G$ on $\mathbb Y_U$ is acylindrical. Thus, applying Theorem \ref{thm:linear_progress_mt} and Proposition \ref{prop:BBF_lowerbound}, we get that:


%
%

\begin{cor}
\label{lem:linear_number_projections} For all $T \geq T_0$, there exist constants $\epsilon_0, C_0>0$ such that: $$ \mathbb P \Big[  \vert \mathcal H_T(1,z_n) \vert  \geq \epsilon_0 n \Big] \geq 1-C_0e^{-n/C_0}. $$
\end{cor}
which is exactly Proposition \ref{prop:Axiom_blackbox}.
%

\subsubsection{Proof of Proposition \ref{prop:axiom_implies_small_interesction} for groups satisfying Assumption \ref{assump:random_walk} }

We first show that the probability that the Gromov product of the position of a random walk $w_n$ at step $n$ taken with respect to a fixed geodesic is larger than $t$ decays exponentially in $t$:

\begin{lemma}
	\label{lem:Gromov_product_decays_in_t}
	Let $G$ be a finitely generated group acting acylindrically on a hyperbolic space $X$, with basepoint $x_0$. Then there exists a constant $D_1>0$ such that the following holds. For all $o \in G$ and $t>0$:
	$$ \mathbb P \Big[(ox_0, w_n x_0)_{x_0} >t \Big] \leq D_1e^{-t/D_1}.$$ 
\end{lemma}

\begin{proof}
	Recall that the shadow $S_{x_0}(x, R)$ is defined (for example in \cite[Section 2]{MaherTiozzo}) as $$ S_{x_0}(x, R)=\{y \in X : (x,y)_{x_0} \geq d_X(x_0,x)-R  \}. $$ 
	
	The following statement is shown in \cite[Lemma 2.10]{MaherMCG} to hold in the mapping class group and in \cite[Equation (16)]{MaherTiozzo} to hold for all acylindrically hyperbolic groups. We will use the formulation of \cite[Proposition 9.3]{MaherSisto}. There 
exists a constant $R_0>0$ only depending on the action of $G$ on $X$, and constants $K>0$ and $D_1'>0$ such that for all $g \in G$ and $R \geq R_0$ : \begin{equation}
\label{eqn:Maher_sisto}
	 \mathbb P\Big[ w_n \in S_{x_0}(gx_0, R) \Big] \leq Ke^{-D_1'(d_X(x_0, gx_0)-R)}.
\end{equation}

Let $D_1 := \max \{Ke^{D_1'R_0}, 1/D_1'\}$ and let $o \in G$ and $t>0$. There are three cases to consider, depending on $o$.

\textbf{Case 1:} If $d_X(x_0, ox_0)-t<0$, then $(ox_0, w_n x_0)_{x_0} <t$, that is,
$$\mathbb P \big[(ox_0, w_n x_0)_{x_0} >t\big]=0.$$

\textbf{Case 2:} If  $d_X(x_0, ox_0)-t \geq R_0$ then by definition of shadows, if $(ox_0, w_n x_0)_{x_0} >t$ then $w_n \in S_{x_0}(ox_0, d_X(x_0, ox_0)-t)=\{y: (ox_0,y)_{x_0}\geq t\}$. Hence, by Equation (\ref{eqn:Maher_sisto}) (as $d_X(x_0, ox_0)-t \geq R_0$): $$\mathbb P \Big[(ox_0, w_n x_0)_{x_0} >t \Big] \leq \mathbb P\Big[w_n \in S_{x_0}(ox_0, d_X(x_0, ox_0)-t) \Big] \leq Ke^{-D_1't} \leq D_1e^{-t/D_1}.$$ 


\textbf{Case 3:} If  $0\leq d_X(x_0, ox_0)-t< R_0$, then $(ox_0, w_n x_0)_{x_0} >t$ implies that $(ox_0, w_n x_0)_{x_0} >d_X(x_0, ox_0)-R_0 $ and so $w_n \in S_{x_0}(ox_0, R_0)$. And hence again by Equation (\ref{eqn:Maher_sisto}) we get:
$$\mathbb P \Big[(ox_0, w_n x_0)_{x_0} >t \Big] \leq \mathbb P\Big[w_n \in S_{x_0}(ox_0, R_0) \Big] \leq Ke^{-D_1'(d_X(x_0, ox_0)-R_0)} \leq Ke^{D_1'R_0}e^{-D_1't} \leq D_1e^{-t/D_1}$$
as $d_X(x_0, ox_0) >t$ and by the choice of $D_1>0$. This finishes the proof of the lemma.
	\end{proof}

We can now prove that Proposition \ref{prop:axiom_implies_small_interesction} holds for random walks satisfying Assumption \ref{assump:random_walk}. 

\begin{proof}[Proof of Proposition \ref{prop:axiom_implies_small_interesction} for groups  satisfying Assumption \ref{assump:random_walk}]
Let $C_1:=\max\{D_1e^{\Delta/D_1},2D_1\}$ where $\Delta$ is from Lemma \ref{lem:gromov_proj} and $D_1$ is from Lemma \ref{lem:Gromov_product_decays_in_t}.

We note that by Lemma \ref{lem:gromov_proj}, if $\vert \mathcal H_T(o,p) \cap \mathcal H_T(p, w_n^p)\vert >t $ then $(ox_0,w_n^{p}x_0)_{px_0}>t/2-\Delta$ in $\mathbb Y_T$, for some basepoint $x_0$. Now, by  
Theorem \ref{thm:action_on_BBF_acylindrical}, the action of $G$ on $\mathbb Y_T$ (with basepoint $x_0$) is acylindrical and hence we can use Lemma \ref{lem:Gromov_product_decays_in_t}: 
$$ \mathbb P \Big[\vert \mathcal H_T(o,p) \cap \mathcal H_T(p, w_n^p)\vert >t  \Big] \leq \mathbb P \Big[ (ox_0,w_n^{p}x_0)_{px_0}>t/2-\Delta\Big] \leq D_1e^{\Delta/D_1}e^{-t/(2D_1)} \leq C_1e^{-t/C_1}$$
by the choice of $C_1$. 
\end{proof}

\subsection{Proof of Proposition \ref{prop:Axiom_blackbox} and Proposition \ref{prop:axiom_implies_small_interesction}  for groups satisfying Assumption \ref{assump:markov_chain}}

In this subsection, we assume that $G$ and the Markov chain $(w_n^{o})_n$ satisfy Assumption \ref{assump:markov_chain}. In particular they satisfy the conclusion of \cite[Proposition 5.1]{GS21}. This proposition states, roughly, that at any given point $p$, the probability that the Markov chain starting at $p$ undoes more than '$t$' worth of projection on these cosets decays exponentially in $t$. We recall the conclusion of the Proposition for the convenience of the reader:
\begin{itemize}
	\item For each $T\geq T'_0$ there exists a constant $C_1'$ such that for all $o,p\in G$,  $n \in \mathbb{N}$, and  $t>0$ we have:
$$ \mathbb{P}\left[\exists r\leq n : \sum_{\mathcal{H}_T(o,p)} [p,w^p_r] \geq t \right] \leq C_1'e^{-t/C_1'}.$$
\end{itemize}

 Fix $T_0 \geq \max\{100(B+\theta+L, T_0'\}$ where $\theta, B, L$ are from Notation \ref{not:all_constants_fixed}. We first prove Proposition \ref{prop:axiom_implies_small_interesction}, as this is an easy consequence of the conclusion of \cite[Proposition 5.1]{GS21}.
 
 \begin{proof}[Proof of Proposition \ref{prop:axiom_implies_small_interesction} for groups satisfying Assumption \ref{assump:markov_chain}]
	We note that if $\vert  \mathcal H_T(x,g) \cap \mathcal H_T(g,w_k^g) \vert > t $ then:
	
	$$\sum_{\mathcal H_T(x,g)}[g, w_k^g] \geq  \sum_{\mathcal H_T(x,g) \cap \mathcal H_T(g,w_k^g)}[g, w_k^g] \geq T\cdot\vert  \mathcal H_T(x,g) \cap \mathcal H_T(g,w_k^g) \vert > Tt.$$ 
	
	Hence by the conclusion of Proposition \cite[Proposition 5.1]{GS21} this leads to:
	$$\mathbb P \Big[\exists k \leq m :\vert  \mathcal H_T(x,g) \cap \mathcal H_T(g,w_k^g) \vert > t  \Big] \leq \mathbb P \Big[\exists k \leq m:\sum_{\mathcal H_T(x,g)}[g, w_k^g] \geq tT \Big] \leq C_1'e^{-tT/C_1'} \leq C_1e^{-t/C_1}  $$
for some constant $C_1>0$. This proves Proposition \ref{prop:axiom_implies_small_interesction}.
\end{proof}

In order to prove Proposition \ref{prop:Axiom_blackbox}, we will actually show something slightly stronger, namely that there exist constants $\epsilon_0$ and $C_0>0$ such that for all $o \in G$, we have:

$$\mathbb P \Big[  \vert \mathcal H_T(o,w^p_n) \vert -\vert \mathcal H_T(o,p) \vert \geq \epsilon_0 n \Big] \geq 1-C_0e^{-n/C_0}.$$

By the coarse triangular inequality (Corollary \ref{cor:triangle_inequality}), this is enough to prove Proposition \ref{prop:Axiom_blackbox}. We first establish the following result, which will allow to do a similar calculation as in \cite[Section 6]{GS21}.

\begin{prop}
\label{prop:small_expectation}
	For all $T \geq T_0$, there exists $\lambda, \kappa>0$ and $m$ such that for all $o,p \in G$:
	$$
\mathbb{E}\Big[\exp\left(\lambda \left(\vert \mathcal H_T(o,p)\vert  -  \vert \mathcal H_T(o,w_m^p)\vert \right)\right)\Big] \leq 1-\kappa.
$$

\end{prop}

We now fix some $T \geq T_0$. The following lemma tells us that the Markov chain is creating some logarithmic number of 'new' cosets.
	
\begin{lemma}
\label{lem:log_number_of_cosets}
There exist $\eta>0$ and a function $h: \mathbb N \to \mathbb N$ with $h(m) \to 0$, such that for all $m$: $$
	\mathbb P \Big[ \exists k \leq m :\vert \mathcal H_T(p,w_k^p) \backslash \mathcal H_T(o,p) \vert \geq \eta \log(m)  \Big] \geq 1-h(m). $$
\end{lemma}

\begin{proof}
	The following claim is similar to \cite[Lemma 3.17]{GS21}. 
	\begin{claim} There exist constants $U, \eta >0$ such that for all $p$:
		$$\mathbb P \Big[ \exists k \leq m : \vert \mathcal H_T(p, w_k^p) \vert \geq \eta \log(m) \Big] \geq 1-Ue^{-\sqrt{m}/U}.$$ 
	\end{claim}
	
	\begin{proof}[Proof of Claim]
	
	Let $\eta, U>0$ be as in \cite[Lemma 2.10]{GS21}, where we note that the proof there does not use the 'non-amenability' criterion from tameness
	 hence it holds for the Markov chains considered here. Also note that the final '$\eta$' in the claim will actually be different to this one. Fix $h$ to be an element of $G$ that is not in $\gamma=E(g)$ and let $\eta ' = \eta /(T+\ell_G(h))$. Let $y:= g^Thg^T\cdots g^Th$ where the number of $g^T$ is $\lceil \eta' \log(m) \rceil$. Then we have that $ \ell_G(y) \leq (\eta' \log(m)+1)(\ell_G(h)+T) \leq \eta \log(m)$ for all $m \geq m_0$ for some $m_0$. Hence, by \cite[Lemma 2.10]{GS21} for all $m \geq m_0$, we have:
	
	$$ \mathbb P \Big[ \exists i,j \leq m : (w^{p}_i)^{-1}w_j^p=g^Thg^T\cdots g^Th\Big] \geq 1-e^{-\sqrt{m}/U}. $$
	
	Hence $\mathbb P \Big[\vert \mathcal H_T(w^p_i, w_j^p) \vert \geq \eta' \log(m)\Big]  \geq 1-e^{-\sqrt{m}/U} $ and so, by the coarse triangular inequality for $\mathcal H_T$, we get that $\mathbb P\Big[\vert \mathcal H_T(p, w_k^p) \vert \geq \eta' \log(m)/2 \Big] \geq 1-e^{-\sqrt{m}/U}$ for either $k=i$ or $k=j$.
	
	By increasing $U$ we can cover the cases where $m <m_0$, and we get that for all $m$: $$\mathbb P \Big[ \exists k \leq m : \vert \mathcal H_T(p, w_k^p) \vert \geq \eta \log(m) \Big] \geq 1-Ue^{-\sqrt{m}/U}$$ proving the claim.
	\end{proof}

 Let $\eta >0$ be as in the claim, we note that the final '$\eta$' in the statement of  Lemma \ref{lem:log_number_of_cosets} will actually be half of this one. We note that by Proposition \ref{prop:axiom_implies_small_interesction} for all $q \in \mathbb N$, we have: $$ \mathbb P[ \vert \mathcal H_T(o,p) \cap \mathcal H_T(p,w_q^p) \vert > \eta \log(m)/2 ] \leq C_1m^{-\eta/2C_1}.$$
	
Let $\mathcal D_{q,g,m}$  be the event "$w_q^p=g$ and  $\vert \mathcal H_T(p,g) \vert \geq \eta \log(m)$ and $\forall i <q: \vert \mathcal H_T(p,w_i^p) \vert < \eta \log(m)$ " (i.e. $w_q^p=g$ and $q$ is the first time such that $\vert \mathcal H_T(p,w_q^p) \vert \geq \eta \log(m)$). Conditioning on this event, for a fixed $m$, we get:

\begin{align*}
\begin{split}
	\mathbb P \Big[ \exists k \leq m :\vert \mathcal H_T(p,w_k^p) &\backslash \mathcal H_T(o,p) \vert \geq \frac{\eta}{2} \log(m)  \Big] =\sum_{\substack{q \leq m \\ g \in G}} \mathbb P \Big[ \exists k \leq m :\vert \mathcal H_T(p,w_k^p) \backslash \mathcal H_T(o,p) \vert \geq \frac{\eta}{2} \log(m)\Big\vert \mathcal D_{q,g,m}\Big]\mathbb P \big[\mathcal D_{q,g,m} \big]  \\
	&\geq \sum_{\substack{q \leq m \\ g \in G}} \mathbb P \Big[\vert \mathcal H_T(p,w_q^p) \cap \mathcal H_T(o,p) \vert \leq \frac{\eta}{2} \log(m) \Big] \mathbb P \big[ \mathcal D_{q,g,m}\big] \\
	&\geq (1-C_1m^{-\eta/2C_1})\sum_{\substack{q \leq m \\ g \in G}}\mathbb P \big[\mathcal D_{q,g,m} \big]  \\
	&= (1-C_1m^{-\eta/2C_1}) \mathbb P \Big[ \exists k \leq m : \vert \mathcal H_T(p, w_k^p) \vert \geq \eta \log(m) \Big] \\
	&\geq (1-C_1m^{-\eta/2C_1})(1-e^{-\sqrt{m}/U})
\end{split}
\end{align*}

where we use Proposition \ref{prop:axiom_implies_small_interesction} to go from the second to the third line and the Claim above to go from the fourth to the fifth line.
	
Therefore: $$ \mathbb P \Big[ \exists k \leq m :\vert \mathcal H_T(p,w_k^p) \backslash \mathcal H_T(o,p) \vert \geq \frac{\eta}{2} \log(m)  \Big] \geq 1-h(m)$$ for some function $h(m) \to 0$. \end{proof}

\begin{lemma}
\label{lem:complement_big}
Let $\eta$ be as in Lemma \ref{lem:log_number_of_cosets}. Then there exists a function $u:\mathbb N \to \mathbb N$ with $u(m) \to 0$, such that: $$
	\mathbb P \Big[ \vert \mathcal H_T(o,w_m^p) \backslash \mathcal H_T(o,p) \vert \geq \eta \log(m)\Big] \geq 1-u(m).$$
\end{lemma}
Before proving Lemma \ref{lem:complement_big}, we need some preliminary results. For all $m$ and $k \leq m $ and $g \in G$, we let $\mathcal E_{k,g,m}$ is the event  $"w_k^p =g$ and  $k$ is the first time such that $\vert \mathcal H_T(p,w_k^p) \backslash \mathcal H_T(o,p)\vert \geq \eta \log(m)"$, we note that for a fixed $m$ and different $k$ and $g$, the events $\mathcal E_{k,g,m}$ are disjoint.
 
\begin{lemma}
\label{lem:small_interesction_implies_big_number_ofcosets} 
	There exist constants $M_0, C_2>0$ such that for all $m \geq M_0$ the following holds. For all $ k \leq m$ and $g \in G$ we have: $$ \mathbb P \Big[ \vert \mathcal H_T(o,w_m^p) \backslash \mathcal H_T(o,p) \vert \geq \eta \log(m)/2 \Big\vert \mathcal E_{k,g,m} \Big]   \geq 1-C_2m^{-\eta/C_2}.$$
\end{lemma}

\begin{proof}
Let $\mathcal S_1, \mathcal S_2 \subseteq \mathcal H_T(p,w_k^p)$ be defined by $\mathcal S_1:= \{h\gamma:\pi_{h\gamma}(w_m^p) =\pi_{h\gamma}(p)\} \cap \mathcal H_T(p,w_k^p)$ and $\mathcal S_2:= \{h\gamma:\pi_{h\gamma}(w_m^p) = \pi_{h\gamma}(w_k^p)\} \cap \mathcal H_T(p,w_k^p)$. Note that $\mathcal S_2 \subseteq \mathcal H_T(p,w_m^p)$.

We note that by Lemma \ref{lem:2cosetsmax} we have that $ \vert \mathcal H_T(p,w_k^p) \vert \leq \vert \mathcal S_1 \vert + \vert \mathcal S_2 \vert +2 $. Hence: $$\vert \mathcal H_T(p,w_m^p)\backslash \mathcal H_T(o,p) \vert \geq \vert \mathcal S_2\backslash \mathcal H_T(o,p) \vert \geq \vert \mathcal H_T(p,w_k^p)\backslash \mathcal H_T(o,p) \vert -\vert \mathcal S_1 \vert -2. $$

Hence by Lemma \ref{lem:bound_number_cosets} we have that:

\begin{align*}
\begin{split}
	\vert \mathcal H_T(o,w_m^p) \backslash \mathcal H_T(o,p) \vert & \geq \vert \mathcal H_T(p,w_m^p)\backslash \mathcal H_T(o,p) \vert -2 \\
	&\geq  \vert \mathcal H_T(p,w_k^p)\backslash \mathcal H_T(o,p) \vert -\vert \mathcal S_1 \vert -4
	\end{split}
\end{align*}

Hence:
\begin{align*}
\begin{split}
	\mathbb P \Big[ \vert \mathcal H_T(o,w_m^p) \backslash \mathcal H_T(o,p) \vert \geq \eta \log(m)/2 \Big\vert \mathcal E_{k,g,m} \Big] &\geq \mathbb P \Big[\vert \mathcal H_T(p,w_{k}^p) \backslash \mathcal H_T(o,p)  \vert -\vert \mathcal S_1\vert-4 \geq \eta \log(m)/2 \vert \mathcal E_{k,g,m} \Big] \\
	&\geq \mathbb P \Big[ \vert \mathcal S_1 \vert  \leq \eta \log(m)/2-4 \vert \mathcal E_{k,g,m}\Big].\\
\end{split}	
\end{align*}
Now, noting that once we have conditioned on $\mathcal E_{k,g,m}$ then $w_k^p=g$ is fixed. Hence by the conclusion of Proposition \cite[Proposition 5.1]{GS21}:
$$ \mathbb P \Big[ \vert \mathcal S_1 \vert  > \eta \log(m)/2-4 \vert \mathcal E_{kg}\Big] \leq \mathbb P \Big[ \sum_{\mathcal H_T(p,g)}[g,w_{m-k}^g]  > T(\eta \log(m)/2-4)\Big] \leq C_1'm^{-T(\eta \log(m)/2-4)} \leq C_2m^{-\eta/C_2}$$
for some constant $C_2>0.$ \end{proof}

We may now use the above lemma to condition on the event $\mathcal E_{k,g,m}$ as follows:

\begin{proof}[Proof of Lemma \ref{lem:complement_big}]

Using the event $\mathcal E_{k,g,m}$ defined above:
	\begin{align*}
 \begin{split}
 \mathbb P \Big[ \vert &\mathcal H_T(o,w_m^p) \backslash \mathcal H_T(o,p) \vert \geq \eta \log(m)/2\Big] =\sum_{\substack{k \leq m\\ g \in G}}  \mathbb P \Big[ \vert \mathcal H_T(o,w_m^p) \backslash \mathcal H_T(o,p) \vert \geq \eta \log(m)/2 \Big\vert \mathcal E_{k,g,m} \Big] \mathbb P\Big[\mathcal E_{k,g,m}\Big] \\
 & \geq (1-C_2m^{-\eta/C_2})\mathbb P \Big[ \exists k \leq m :\vert \mathcal H_T(p,w_k^p) \backslash \mathcal H_T(o,p) \vert \geq \eta \log(m)  \Big] \\
 & \geq (1-C_2m^{-\eta/C_2})(1-h(m))=1-u(m)
 \end{split}
 \end{align*}
 where we go from the second to the third line by using Lemma \ref{lem:log_number_of_cosets} and Lemma \ref{lem:small_interesction_implies_big_number_ofcosets}, for some function $u$ with $u(m) \to 0$.
\end{proof}

\begin{proof}[Proof of Proposition \ref{prop:small_expectation}]

By the Cauchy-Schwarz inequality, we get:
	$$
		\mathbb{E}\left[\exp\Big(\lambda \left(\vert  \mathcal{H}_T(o,p) \vert -\vert \mathcal H(o,w^p_m) \vert\right)\Big)\right] \leq \mathbb{E}\left[\exp\left(2\lambda \left(\vert  \mathcal{H}_T(o,p)\backslash \mathcal H_T(o,w_m^p) \vert \right)\right)\right]^{1/2}\mathbb{E}\left[\exp\left(-2\lambda \left(\vert  \mathcal{H}_T(o,w_m^p)\backslash \mathcal H_T(o,p) \vert \right)\right)\right]^{1/2}. 
$$

We calculate each of these terms individually.
The first expected value can be computed using the following claim.

\begin{claim}
\label{claim:small_proba_complementH(o,p)}
	There exists a constant $C_3>0$ such that for all $o,p \in G$, $m \in \mathbb N$ and $s\geq 0$ we have:
	
$$\mathbb P\Big[ \vert \mathcal H_T(o,p)\backslash \mathcal H_T(o,w_m^p) \vert \geq s \Big] \leq C_3e^{-s/C_3}. $$
\end{claim}

\begin{proof}
By Lemma \ref{lem:2cosetsmax} we have at most 2 cosets (call them $t_1\gamma, t_2\gamma$) from $\mathcal H_T(o,p)\backslash \mathcal H_T(o,w_m^p)$ such that $\pi_{t_i\gamma}(w_m^p) \neq \pi_{h\gamma}(p)$ and $\pi_{t_i\gamma} \neq \pi_{h\gamma}(o)$. For all but these 2 cosets $t_i \gamma$, we must have $\pi_{h\gamma}(w_m^p)=\pi_{h\gamma}(o)$ (otherwise $h\gamma \in  \mathcal H_T(o,w_m^p)$). Therefore, for all cosets  $h\gamma \in \mathcal H_T(o,p)\backslash \mathcal H_T(o,w_m^p)$ (apart from $t_i\gamma$) we have $d_{h\gamma}(p,w_m^p)=d_{h\gamma}(p,o) \geq T$.  Hence if 
$\vert \mathcal H_T(o,p)\backslash \mathcal H_T(o,w_m^p) \vert \geq s$ then 

$$\sum_{\mathcal H_T(o,p)}[p,w_m^p] \geq \sum_{\mathcal H_T(o,p)\backslash \mathcal H_T(o,w_m^p)}d_{h\gamma}(p,w_m^p) \geq \left(\vert \mathcal H_T(o,p)\backslash \mathcal H_T(o,w_m^p) \vert-2\right)T \geq (s-2)T $$

and so by the conclusion of \cite[Proposition 5.1]{GS21}: $$ \mathbb P\Big[ \vert \mathcal H_T(o,p)\backslash \mathcal H_T(o,w_m^p) \vert \geq s\Big]\leq \mathbb P\Big[ \sum_{\mathcal H_T(o,p)}[p,w_m^p] \geq (s-2)T\Big] \leq C_1'e^{-(s-2)T/C_1'}.$$

Increasing the constant to a constant $C_3$ such that $C_1'e^{-(s-2)T/C_1'}\leq C_3e^{-s/C_3}$ proves the Claim.
\end{proof}

Hence, using Claim \ref{claim:small_proba_complementH(o,p)} and choosing $\lambda<1/(2C)$, we can bound the first term in the expected value above as follows:
\begin{align*}
	\begin{split}
		\mathbb{E}\Bigg[\exp\Big(2\lambda  \big(\vert  &\mathcal{H}_T(o,p)\backslash \mathcal H_T(o,w_m^p) \vert \big)\Big)\Bigg]
=\int_{0}^{+\infty} \mathbb P \Big[ \exp(2\lambda \left(\vert  \mathcal{H}_T(o,p)\backslash \mathcal H_T(o,w_m^p) \vert \right) \geq s\Big]ds \\
		&=\int_{0}^{+\infty} \mathbb P\Big[ \vert  \mathcal{H}_T(o,p)\backslash \mathcal H_T(o,w_m^p) \vert \geq \ln(s)/2\lambda\Big] \\
		&=\int_{0}^{1} \mathbb P\Big[ \vert  \mathcal{H}_T(o,p)\backslash \mathcal H_T(o,w_m^p) \vert \geq \ln(s)/2\lambda\Big]+\int_{1}^{+\infty} \mathbb P\Big[ \vert  \mathcal{H}_T(o,p)\backslash \mathcal H_T(o,w_m^p) \vert \geq \ln(s)/2\lambda\Big]\\
		&\leq 1+C\int_{1}^{+\infty} s^{-1/(2\lambda C)} ds \\
		&=1+\frac{2C^2\lambda}{1-2\lambda C}.
	\end{split}
\end{align*}
This bounds the first term in the expression for the expected value.
We can find an upper bound for the second term  as follows. Let $A$ be the event: "$\vert  \mathcal{H}(o,w_m^p)\backslash \mathcal H(o,p) \vert \geq \eta \log(m)/2$", then by Lemma \ref{lem:complement_big}, we get

$$\mathbb{E}\Big[\exp\Big(-2\lambda \left(\vert  \mathcal{H}_T(o,w_m^p)\backslash \mathcal H_T(o,p) \vert \right)\Big)\Big]= $$
\begin{align*}
\begin{split}
	\mathbb{E}\left[\exp\Big(-2\lambda \left(\vert  \mathcal{H}_T(o,w_m^p)\backslash \mathcal H_T(o,p) \vert \right)\mathds{1}_A\right]+\mathbb{E}\left[\exp\Big(-2\lambda \left(\vert  \mathcal{H}_T(o,w_m^pp)\backslash \mathcal H_T(o,p) \vert \right)\mathds1_{A^C}\right]  &\leq m^{-\lambda\eta}+ \mathbb P [ A^C]\\
	&\leq m^{-\lambda\eta}+u(m).\\
\end{split}	
\end{align*}
Therefore:
\begin{align*}
\begin{split}
	\mathbb{E}\left[\exp\Big(\lambda \left(\vert  \mathcal{H}_T(o,p) \vert -\vert \mathcal H_T(o,w^p_m) \vert\right)\Big)\right] \leq \left(1+\frac{2C^2\lambda}{1-2\lambda C}\right)^{1/2}\left(u(m)+m^{-\lambda\eta}\right)^{1/2}
\end{split}	
\end{align*}
Recall that $\lambda$ has been fixed, there is $m$ and $\kappa >0$ such that: 
$$ \mathbb{E}\left[\exp\Big(\lambda \left(\vert  \mathcal{H}_T(o,p) \vert -\vert \mathcal H_T(o,w^p_m) \vert\right)\Big)\right] \leq 1-\kappa.$$\end{proof}
We now have all the necessary results in order to prove Proposition \ref{prop:Axiom_blackbox} for groups and Markov chains satisfying Assumption \ref{assump:markov_chain}.
\begin{proof}[Proof of Proposition \ref{prop:Axiom_blackbox} for groups satisfying Assumption \ref{assump:markov_chain}]

We have for all integers $j>0$: 

\begin{align*}
\begin{split}
	 \mathbb{E}\Big[\exp\Big(\lambda &(\vert  \mathcal{H}_T(o,p) \vert -\vert \mathcal H_T(o,w^p_{(j+1)m}) \vert)\Big)\Big] \\
  &=\mathbb{E}\Big[\exp\Big(\lambda (\vert  \mathcal{H}_T(o,p) \vert -\vert \mathcal H_T(o,w^p_{(j+1)m}) \vert-\vert\mathcal H_T(o,w^p_{jm}) \vert+\vert \mathcal H_T(o,w^p_{jm})\vert)\Big)\Big] \\
	 &= \sum_{g \in G} \mathbb{E}\Big[e^{(\lambda ( \vert  \mathcal{H}_T(o,p) \vert -\vert \mathcal H_T(o,w^p_{(j+1)m}) \vert-\vert \mathcal H_T(o,w^p_{jm}) \vert+\vert \mathcal H_T(o,w^p_{jm})\vert)} \Big\vert w_{jm}=g\Big]\mathbb P \Big[w_{jm}=g \Big]\\
	 &\leq (1-\kappa) \mathbb E \Big[ e^{-\lambda(\vert \mathcal H_T(o,w^p_{jm})\vert-\vert \mathcal{H}_T(o,p) \vert)} \Big] \\
\end{split}
\end{align*}

Hence, by induction, we get that 
$$ \mathbb{E}\left[\exp\Big(\lambda \left(\vert  \mathcal{H}_T(o,p) \vert -\vert \mathcal H_T(o,w^p_{jm}) \vert\right)\Big)\right] \leq (1-\kappa)^j. $$

Therefore, by the Markov inequality:

\begin{align*}
\begin{split}
 \mathbb P \Big[\vert \mathcal H_T(o,w_{jm}^p) \vert - \vert \mathcal H_T(o,p) \vert \leq \epsilon_0 jm \Big] &=\mathbb P \Big[\exp\left(-\lambda(\vert \mathcal H_T(o,w_{jm}^p) \vert - \vert \mathcal H_T(o,p) \vert \right)) >\exp(-\lambda \epsilon_0 jm) \Big]  \\
 &\leq \frac{\mathbb{E}\left[\exp\Big(\lambda \left(\vert  \mathcal{H}_T(o,p) \vert -\vert \mathcal H_T(o,w^p_{jm}) \vert\right)\Big)\right]}{\exp(-\lambda \epsilon jm)} \\
 &\leq (1-\kappa)^{j}e^{\lambda \epsilon_0 jm}.
\end{split}
\end{align*}

Choosing $\epsilon_0$ small enough, we can find a constant $C_0$ such that 

$$ \mathbb P \Big[\vert \mathcal H_T(o,w_{jm}^p) \vert - \vert \mathcal H_T(o,p) \vert \leq \epsilon_0 jm \Big] \leq C_0e^{-jm/C_0}. $$

By a similar discussion to the one at the start of  \cite[Section 6]{GS21}, this is enough to prove that for all $n$, we have:
$$\mathbb P \Big[ \vert \mathcal H_T(o,w_n^p) \vert -\vert \mathcal H_T(o,p) \vert \geq \epsilon_0 n \Big] \geq 1-C_0e^{-n/C_0}. $$

As we mentioned earlier, by the coarse triangular inequality (Corollary \ref{cor:triangle_inequality}), this proves Proposition \ref{prop:Axiom_blackbox}. 

Hence, we have established that both Proposition \ref{prop:Axiom_blackbox} and \ref{prop:axiom_implies_small_interesction} hold if the group $G$ and the Markov chain $(w^{o}_n)_n$ satisfy Assumption \ref{assump:markov_chain}.
\end{proof}
\section{Proof of Theorem \ref{thm:lower_bound}}
\label{sec:main_proof}

We now prove that if Proposition \ref{prop:Axiom_blackbox} and Proposition \ref{prop:axiom_implies_small_interesction} hold then Theorem \ref{thm:lower_bound} holds. Fix $T$ be greater than both $T_0$ and $T_0'$ from Proposition \ref{prop:Axiom_blackbox} and Proposition \ref{prop:axiom_implies_small_interesction}. Let $\epsilon_0, C_0>0$  be as in Proposition \ref{prop:Axiom_blackbox} for this $T$.  For all $x,y \in G$ and $S >0$, we let $$\mathcal A^{S}_{x,y} := \mathcal H_T (x,y) \cap \{ h \gamma : \pi_{h\gamma}(x) \subseteq B(x, S) \}.  $$

\begin{remark}
	We note that $\mathcal H_T(x,y) =\mathcal H_T(y,x)$ but in general $\mathcal A_{x,y}^S \neq \mathcal A_{y,x}^S $.
\end{remark}

The following proposition is crucial when proving Theorem \ref{thm:lower_bound}, and states that not only does the number of cosets in $\mathcal H_T(p,w_n^p)$ grow linearly but also those cosets where $\pi_{h\gamma}(p)$ is close to $p$.

\begin{prop}
\label{prop:linear_implies_in_ball}
For all $\epsilon \leq \epsilon_0$,	there exist constants $\nu, C_4 >0$ such that for all $n \in \mathbb N$: 
	
 $$ \mathbb P \Big[ \vert \mathcal A^{\epsilon n/2}_{p,w^p_n} \vert \geq \nu n  \Big] \geq 1-C_4e^{-n/C_4}.$$ 
\end{prop}

In order to prove Proposition \ref{prop:linear_implies_in_ball}, we need various lemmas. We let $\epsilon \leq \epsilon_0$ which we fix throughout this proof. Then by Proposition \ref{prop:Axiom_blackbox}:
\begin{equation}
\label{eqn:linear_number_cosets}
\mathbb P \Big[ \vert \mathcal H_T(p,w^{p}_n) \vert \geq \epsilon n \Big]  \geq \mathbb P \Big[  \vert \mathcal H_T(p,w^{p}_n) \vert \geq \epsilon_0 n \Big] \geq 1-C_0e^{-n/C_0}. \end{equation}
We now fix $n \in \mathbb N$. For ease of notation, as $\epsilon$ has been fixed, we will write $\mathcal A_{x,y} := \mathcal A^{\epsilon n/2}_{x,y}$, for all $x,y \in G$. For any $U \geq 0$, we let $\mathcal B_{U} \subseteq G$ be the set of elements such that:
$$ \mathcal B_U := \{g \in G : \vert \mathcal A_{p,g} \vert \geq U  \}. $$

For all $k \leq n$ and $g \in G$, we consider the event $B_{k, U,g}:" w_k^p=g$ and $g \in \mathcal B_U$ and for all $i <k$ : $w_i^p \not\in \mathcal B_U$" (i.e. $w_k^p=g$ is the first time $k$ that $w_k^p \in B_U$).  We note that for any $U \geq 0$ and for distinct  $k_1, k_2 \in \mathbb N $,  the events $ B_{k_1, U,g}$ and $B_{k_2, U,g}$  are disjoint. We need the following two lemmas.
\begin{lemma}
	\label{lem:exists_initial_goodpath}
	
There exist constants $\epsilon_3, C_5 >0$ such that $$  \mathbb P \Big[ \exists k \leq n :  \vert \mathcal A_{p,w^p_k} \vert \geq \epsilon_3 n \Big] \geq 1-C_5e^{-n/C_5}. $$
\end{lemma}

The following lemma will later allow us to condition on the event defined above.

\begin{lemma}
	\label{lem:consequence_of_small_gromov}
	For all $\tau>0$ there exists a constant $C_6$ such that the following holds. For all $k \leq n$, we have:
	
$$\mathbb P \Big[ \vert \mathcal A_{1,z_n} \vert \geq \frac{9\tau n}{10} \quad  \Big\vert B_{k,\tau n,g}\Big] \geq 1-C_6e^{-n/C_6}. $$
\end{lemma}

It now remains to prove these two lemmas.  We will use the following lemma twice in order to prove the two other ones, so we state it now.

Recall that $\theta, B, K, L$ and $f$ have been fixed in Notation \ref{not:all_constants_fixed}. 
We note that the following lemma is about discrete paths $\alpha$ and  that the neighbourhood $\mathcal N^G_s$ in the lemma below is taken with respect to the distance $d_G$. 

\begin{lemma}
\label{lem:argument_divergence_path_ouside_ball}
	Let $x,y \in G$, let $\alpha$ be a discrete path from $x$ to $y$, let $r>0$, and let $\mathcal F \subseteq \mathcal H_T(x,y)$ be such that for all $h\gamma \in \mathcal F$ the following holds: 
	\begin{itemize}
		\item if $p \in \alpha$ and $d_{h\gamma}(x,p) \leq \theta+B+2L$ then $p \not\in \mathcal N^G_{s}(h\gamma)$.
	\end{itemize}
	Then $\ell(\alpha) > \vert \mathcal F \vert f(s).$
\end{lemma}	
\begin{proof}
Let $\alpha, s, \mathcal F$ as in the statement of the Lemma. With a slight abuse we will regard $\alpha$ as a discrete path, meaning a sequence of points where consecutive points are at distance $1$ from each other, which can be arranged by replacing $\alpha$ with a suitable discretisation. We order the cosets in $\mathcal F$ according to the linear order from $\mathcal H_T(x,y)$. Starting from the smallest such coset, for all $h_i\gamma \in \mathcal F $  we define $\alpha^{-}_i$ as the last point along $\alpha$ such that $d_{h_i\gamma}(x, \alpha^{-}_i) \leq B+L$ and $\alpha^{+}_i$ as the first time that $d_{h_i\gamma}(\alpha^{-}_i, \alpha^{+}_i) \geq \theta$. By the choice of $T\gg \theta,B,L$, this is possible. Let $\alpha_i=\alpha\vert_{[\alpha^{-}_i, \alpha^{+}_i]}$. Note that any $p$ on $\alpha_i$ satisfies $d_{h_i\gamma}(p, x) >B$ (otherwise the point $p'$ after $p$ along $\alpha$ would satisfy $d_{h_i\gamma}(x, p') \leq B+L$ contradicting that $\alpha^{-}_i$ is the last such point) and $d_{h_i\gamma}(p, y) >B$ (otherwise the previous point $p'$ would satisfy $d_{h_i\gamma}(\alpha^{-}_i, p') \geq \theta$). In particular, by Lemma \ref{lem:linearorder}, for all $i\neq j$ we have $d_{h_j\gamma}(p,x)\leq B$ or $d_{h_j\gamma}(p,y)\leq B$ for all $p\in \alpha_i$. This shows that the $\alpha_i$ are pairwise disjoint. 
For all $p\in \alpha_i$ we have $d_{h_i\gamma}(x,p) \leq \theta+B+2L$, and therefore we have $\alpha_i \cap \mathcal N_{s}(h_i\gamma) =\emptyset$ by assumption.  Hence, by $f$-divergence, this leads to:
$$\ell(\alpha) \geq \sum_{h_i\gamma \in \mathcal F } \ell (\alpha_i) > \vert \mathcal F\vert  f(s).$$\end{proof}
Before proving Lemma \ref{lem:exists_initial_goodpath}, we state a few more lemmas. The first lemma is similar to  \cite[Lemma 3.6]{GS21} and ensures that we can apply Lemma \ref{lem:argument_divergence_path_ouside_ball} later. The ball $B^G$ is considered with respect to the metric $d_G$.

\begin{lemma}
\label{lem:avoid_neighbourhood}
There exist a constant $D_2 >0$ and $N \in \mathbb N$ such that the following holds. For all $n \geq N$ and for all $x,z \in G$ such that $z \in B^G(x, \epsilon n/3)$, we have:  $$d_G(x,z) > \vert \mathcal H_T(x,z) \cap \{h\gamma:\pi_{h\gamma}(x) \not\subseteq B^G(x, \epsilon n/2)  \} \vert \cdot f(n/D_2). $$
\end{lemma}

\begin{proof}
Let $D_1 \geq \max\{D,a ,b\}$ where $D$ is from Lemma \ref{lem:distances_projection} and $a,b$ are the quasi-isometric embedding constants $h\gamma \hookrightarrow G$. Let $D_2$ be big enough such that $\epsilon/6-D_1/D_2 >0$ where $D$ is from Lemma \ref{lem:distances_projection} and let $N$ be such that for all $n \geq N$ we have that $D_1(\theta+B+2L)+2D_1 \leq n(\epsilon/6-D_1/D_2)-2B$.

Fix $n \geq N$ and $z\in B(x, \epsilon n/3)$. Let $\mathcal F :=\mathcal H_T(x,z) \cap \{h\gamma:\pi_{h\gamma}(x) \not\subseteq B(x, \epsilon n/2)  \}$. By Lemma \ref{lem:argument_divergence_path_ouside_ball}, it suffices to show for all $h\gamma \in \mathcal F$ and all $p \in [x,z]$ where $d_{h\gamma}(x,p) \leq \theta +B+2L$ we have that $p \not\in \mathcal N_{n/D_2}(h\gamma)$.

Hence, let $h\gamma \in \mathcal F$ and $p \in [x,z]$ such that $d_{h\gamma}(x,p) \leq \theta+B+2L$. For a contradiction, assume  that $p \in \mathcal N_{n/D_2}(h\gamma)$. Then, by the choice of $D_2 $ and $n$ and Lemma \ref{lem:distances_projection} we have: \begin{align*}
 \begin{split}
 	d_G(x, \pi_{h\gamma}(x)) &\leq d_G(x, p)+d_G(p, \pi_{h\gamma}(p)) + d_G(\pi_{h\gamma}(p), \pi_{h\gamma}(x)) \\
 	&\leq \epsilon n/3 + nD_1/D_2+D_1+D_1(\theta+B+2L)+D_1+2B \\
 	&< n\Big(\epsilon/3+\epsilon/6+D_1/D_2-D_1/D_2\Big) = \epsilon n/2\\
 \end{split}	
 \end{align*}
  a contradiction with $\pi_{h\gamma}(x) \not\subseteq B(x, \epsilon n	/2)$ (where the "$2B$" compensates for the diameter of $\pi_{h\gamma}(p)$).
\end{proof}

\begin{notation}
\label{not:fixing_epsilons}
Let $\epsilon_1 < \epsilon/3K$ where $K$ is from the bound on $d_G(w^{p}_k, w^p_{k+1})$ for the Markov chain  $(w^{o}_n)$ (see Notation \ref{not:all_constants_fixed} and Assumption \ref{assump:markov_chain}), fix $\epsilon_2 =\epsilon \cdot\epsilon_1$ and $\epsilon_3=\epsilon \cdot\epsilon_1/3$.
\end{notation}

\begin{lemma}
\label{lem:conditioning_on_k}
	There exists $N_0 \in \mathbb N$ such that for all $n \geq N_0$ the following holds. If $d_G(x,y) \leq \epsilon n /3 $ and $\vert \mathcal H_T(x, y) \vert > \epsilon_2 n$ then $\vert \mathcal A_{x, y} \vert \geq \epsilon_3 n$, where $\epsilon_1, \epsilon_2, \epsilon_3$ are as in Notation \ref{not:fixing_epsilons}.
\end{lemma}

\begin{proof}
Let $D_2>0$ and $N \in \mathbb N$ be from Lemma \ref{lem:avoid_neighbourhood}. Let $n_0$ be such that for all $n \geq n_0$ we have $f(n/D_2) >1/(2\epsilon_1)$ (as $f$ is divergent). Let $N_0 = \max\{n_0, N\}$ and $n \geq N_0$. As $d_G(x,y)\leq \epsilon n/3$ and $n \geq N$, we can apply Lemma \ref{lem:avoid_neighbourhood} to get: $$ \epsilon n/3 \geq d_G(x,y) > \vert \mathcal H_T(x,y) \cap \{h\gamma:\pi_{h\gamma}(x) \not\subseteq B(x, \epsilon n/2)  \} \vert \cdot f(n/D_2) $$ and hence $$\vert \mathcal H_T(x,y) \cap \{h\gamma:\pi_{h\gamma}(x) \not\subseteq B(x, \epsilon n/2)  \} \vert < \frac{\epsilon n}{3f(n/D_2)}. $$
We then get:
		\begin{align*}
			\begin{split}
			\vert \mathcal A_{x, y} \vert &= \vert \mathcal H_T(x, y) \vert -\vert \mathcal H_T(x,y) \cap \{h\gamma:\pi_{h\gamma}(x) \not\subseteq B(x, \epsilon n/2)  \} \vert \\
			& \geq \epsilon_2 n-\frac{\epsilon n}{3f(n/D_2)} \\
			&=\epsilon n\Big(\epsilon_1-\frac{1}{3f(n/D_2)}\Big) \\
			&\geq \epsilon \cdot \epsilon_1 n/3 = \epsilon_3 n. \\
			\end{split}
		\end{align*}	
		\end{proof}
We are now ready to prove Lemma \ref{lem:exists_initial_goodpath}. 		
\begin{proof}[Proof of Lemma \ref{lem:exists_initial_goodpath}]
Let $\epsilon_1, \epsilon_2, \epsilon_3$ be as in Notation \ref{not:fixing_epsilons}. We note that if $k \leq \epsilon_1 n$ then by the $K$-bounded jumps of the Markov chain we get that $d_G(p,w^p_k) \leq K\epsilon_1 n \leq \epsilon n/3$. Let $N'_0$ be greater than $N_0$ in  Lemma \ref{lem:conditioning_on_k} and such that $\epsilon_1N'_0+1 \leq N'_0$. Then for $ n \geq N'_0$ we have:
\begin{align*}
\begin{split}
\mathbb P \Big[ \exists k \leq n :  \vert \mathcal A_{p,w^p_k} \vert \geq \epsilon_3 n \Big]  
 & \geq  \mathbb P \Big[\exists k \leq \epsilon_1 n +1 :\vert \mathcal H_T(p, w^p_k) \vert > \epsilon_2 n  \Big] \\
 &\geq \mathbb P \Big[\vert \mathcal H_T(p, w^p_{\lceil \epsilon_1 n\rceil}) \vert > \epsilon_2 n  \Big] \\
 &\geq 1-C_0'e^{-\lceil \epsilon_1 n\rceil/C_0'}
\end{split}
\end{align*}
by Equation (\ref{eqn:linear_number_cosets}). Up to increasing the constant to some $C_5$ to cover the case where $n \leq N'_0$, this completes the proof of Lemma \ref{lem:exists_initial_goodpath}.
\end{proof}

We now prove the other crucial lemma.
\begin{proof}[Proof of Lemma \ref{lem:consequence_of_small_gromov}]

Note that, for all $g$, we have $\mathcal A_{p,w^{g}_{n-k}} := \mathcal H_T (p,w^{q}_{n-k}) \cap \{ h \gamma : \pi_{h\gamma}(p) \subseteq B(p, \epsilon n/2) \}  $, hence if $h\gamma  \in \mathcal A_{p,w^{p}_k}$  and $h\gamma \in \mathcal H_T(p, w^{g}_{n-k})$ then $h\gamma \in \mathcal A_{p, w^{g}_{n-k}}$.

Therefore for all $g \in G$, $k \leq n$ and $\tau >0$:
\begin{align*}
\begin{split}
	\mathbb P \Big[ \vert \mathcal A_{p,w_{n}^p} \vert \geq 9\tau n /10 \quad  \Big\vert \quad \big( \vert \mathcal A_{p,w^p_k}\vert \geq \tau n\big) \cap (w_k^p=g) \Big] &\geq \mathbb P \Big[ \vert \mathcal H_T(g, w_{n-k}^g) \cap \mathcal H_T(g,p) \vert  \leq \tau n/10\Big] \\
	& \geq 1-C_1e^{-\tau n/(10C_1)} =1-C_6e^{-n/C_6}
\end{split}	
\end{align*}

where we use the strong Markov property (see \cite[Lemma 2.2]{GS21}) and $C_1$ is from Proposition \ref{prop:axiom_implies_small_interesction}, the result follows for some constant $C_6>0$.
\end{proof}

Recall that we have defined the set $\mathcal B_U$ and the event $B_{k,U,g}$ above.
		 
\begin{proof}[Proof of Proposition \ref{prop:linear_implies_in_ball}] 
Recall that we have fixed $\epsilon \leq \epsilon_0$. Let $\epsilon_3, C_5>0$ be from Lemma \ref{lem:exists_initial_goodpath}, let $C_6$ be the constant associated to $\nu=\frac{9\epsilon_3}{10}$ from Lemma \ref{lem:consequence_of_small_gromov}. 
We then have: 
\begin{align*}
	\begin{split}
		\mathbb P \Big[ \vert \mathcal A_{p,w^p_n} \vert \geq \frac{9\epsilon_3}{10}n  \Big] 	& \geq \sum_{\substack{ k \leq n \\ g \in G }} \mathbb P \Big[ \vert \mathcal A_{p,w^p_n} \vert \geq \frac{9\epsilon_3}{10}n \Big\vert B_{k,\epsilon_3 n,g}  \Big] \mathbb P\Big[ B_{k,\epsilon_3 n,g} \Big] \\ 
		& \geq \Big(1-C_6e^{-n/C_6}\Big)\sum_{\substack{ k \leq n }} \mathbb P\Big[ B_{k,\epsilon_3 n,g} \Big] \\ 
		&= \Big(1-C_6e^{-n/C_6}\Big)\mathbb P \Big[ \exists k \leq n :  \vert \mathcal A_{p,w^p_k} \vert \geq \epsilon_3 n \Big] \\
		& \geq \Big(1-C_6e^{-n/C_6}\Big)(1-C_5e^{-n/C_5}) \\
		& \geq 1-C_4e^{-n/C_4}
	\end{split}
\end{align*}

where we go from the first line to the second line by using Lemma \ref{lem:consequence_of_small_gromov}, from the third line to the fourth by using Lemma \ref{lem:exists_initial_goodpath}, and where $C_4$ is sufficiently large depending on $C_5,C_6$.
\end{proof}

Let $D>0$ be from Lemma \ref{lem:distances_projection}. The following deterministic lemma will allow us to find the length of a path avoiding a ball, and hence the divergence as in the Theorem \ref{thm:lower_bound}.
\begin{lemma}
	\label{lem:implies_long_length}
	
For all $\delta>0$, there exists $N_{\delta} \in \mathbb N$ such that for all $r \geq N_{\delta}$ and for all $x,w,z \in G$  the following holds. Any discrete path $\beta$ from $w)$ to $z$ avoiding the ball $B(x, \delta r)$ has length $$\ell (\beta) > \left(\vert \mathcal A^{ \delta r/2}_{x,z} \backslash \mathcal A^{ \delta r/2}_{x,w} \vert -2 \right) \cdot f\big(\frac{\delta r}{4D}\big).$$
	
	\end{lemma}
	
		
	\begin{proof}
		Let  $N'_{\delta}$ be such that for all $r \geq N'_{\delta}$, we have $\delta r /2  \geq \delta r/3+\theta+3B+2L$ and let $N_{\delta} := \max\{N'_{\delta}, 12D/\delta \}$. 
	Fix $n \geq N_{\delta}$. By Lemma \ref{lem:bound_number_cosets} there are at most 2 cosets $t_1\gamma$ and $ t_2\gamma$ in $\mathcal H_T(x,z)$  such that $$  \mathcal H_T(x,z) \backslash \mathcal H_T(x,w) \subseteq \mathcal H_T(z,w) \cup \{ t_1\gamma, t_2\gamma\}.$$ 
 Now, if $h\gamma \in \mathcal A^{\delta r/2}_{x,z} \backslash \mathcal A^{\delta r/2}_{x,w}$ then $\pi_{h\gamma}(x) \subseteq B(x,\delta r/2)$ and hence if $h\gamma \not\in \mathcal A^{\delta r/2}_{x,w}$, this is because $h\gamma \not\in \mathcal H_T(x,w)$. We also clearly have that $\mathcal H_T(x,w)^C \subseteq (A^{\delta r/2}_{x,w})^C$ and hence $\mathcal A^{\delta r/2}_{x,z} \backslash \mathcal A^{\delta r/2}_{x,w}=\mathcal A^{\delta r/2}_{x,z} \backslash \mathcal H_T(x,w)$. Therefore, letting $\mathcal F := \mathcal A^{\delta r/2}_{x,z} \backslash \Big(\mathcal A^{\delta r/2}_{x,w} \cup \{t_1\gamma, t_2\gamma \} \Big)$, we get: $$ \mathcal F =\left( \mathcal A^{\delta r/2}_{x,z}\backslash \mathcal H_T(x,w)\right) \cap  (\{t_1\gamma, t_2\gamma \})^C \subseteq \left( \mathcal H_T(x,z)\backslash \mathcal H_T(x,w)\right) \cap  (\{t_1\gamma, t_2\gamma \})^C  \subseteq \mathcal H_T(z,w). $$
	
Now that we know that $\mathcal F \subseteq \mathcal H_T(z,w) $, we show that it satisfies the other assumptions of Lemma \ref{lem:argument_divergence_path_ouside_ball} for the path $\beta$ and for $s=\frac{\delta r}{4D}$. We note that for all $h\gamma \in \mathcal F$, we have $\pi_{h\gamma}(x)=\pi_{h\gamma}(w)$. Let $p \in \beta$ and $h\gamma \in \mathcal F$ be such that $d_{h\gamma}(x,p) \leq  
	\theta+B+2L$. Then: 
\begin{align*}
	\begin{split}
		d_G(p, \pi_{h\gamma}(p)) &\geq d_G(x,p)-d_G(x, \pi_{h\gamma}(x)) -d_G(\pi_{h\gamma}(x), \pi_{h\gamma}(p)) \\
		& \geq \delta r-\delta r /2-(\theta+B+2L)-2B \\
		&\geq \delta r /3
	\end{split}
\end{align*}
by the choice of $N'_{\delta}$. By Lemma \ref{lem:distances_projection} and the choice of $N_{\delta}$, this leads to $d_G(p,h\gamma) \geq \delta  r/(3D)-1 \geq \delta r/(4D)$. Therefore, we can apply Lemma \ref{lem:argument_divergence_path_ouside_ball} to get:	$$ \ell (\beta)  > \left(\vert \mathcal A^{\delta r/2}_{x,z} \backslash \mathcal A^{ \delta r/2}_{x,w} \vert -2 \right) \cdot f\big(\frac{\delta r}{4D}\big).$$\end{proof}

Now that we can find a lower bound for any path joining the geodesics $[p, w_n^p]$ to $[p,z_n^p]$ and avoiding a ball around $p$, we just need to find a lower bound on the probability that both $\vert \mathcal A^{\delta n/2}_{p,z^{p}_n} \backslash \mathcal A^{\delta n/2}_{p,w^p_n} \vert$ and $d(p, \{w_n^{p}, z_n^p\}) $ are big (in order to  apply Lemma \ref{lem:argument_divergence_path_ouside_ball}).

\begin{proof}[Proof of Theorem \ref{thm:lower_bound}]


Let $\delta_0=\epsilon_0$ and $\delta \leq \delta_0$, then by Proposition \ref{prop:linear_implies_in_ball} there exist constants $\nu, C_4$ (depending on $\delta$) such that $\mathbb P \Big[ \vert \mathcal A^{\delta  \kappa n/2}_{p,z^p_n} \vert \geq \nu n\Big] \geq 1-C_4e^{-n/C_4}. $ Let $0<c< \min\{\nu/2, \delta_0\}$.
We note that:
\begin{align*}
	\begin{split}
		\mathbb P \Big[\vert \mathcal A^{\delta \kappa n/2}_{p,z^{p}_n} \backslash \mathcal A^{\delta \kappa n/2}_{p,w^p_n} \vert -2 \geq cn   \Big] & =\sum_{h \in G} \mathbb P \Big[\vert \mathcal A^{\delta \kappa n/2}_{p,z^{p}_n} \backslash \mathcal A^{\delta \kappa n/2}_{p,w^p_n} \vert -2 \geq cn \Big\vert w_n^p=h \Big] \mathbb P[w_n^p=h] \\
		& \geq \sum_{h \in G} \mathbb P \Big[ \left(\vert \mathcal A^{\delta \kappa n/2}_{p,z_n^p} \vert \geq \nu n \right) \cap \left(\vert \mathcal H_T(p,z_{n}^p) \cap \mathcal H_T(p,h) \vert \leq \nu n/2 \right)\Big] \mathbb P[w_n^p=h] \\
		&\geq \sum_{h\in G} \left(  1-\mathbb P \Big[ \vert \mathcal A^{\delta \kappa n/2}_{p,z_n^p} \vert <\nu n \Big] - \mathbb P \Big[ \vert \mathcal H_T(p,z_{n}^p) \cap \mathcal H_T(p,h) \vert > \nu n/2\Big] \right) \mathbb P[w_n^p=h] \\
		&\geq 1-C_4e^{-n/C_4}-C_1e^{-\nu n/(2C_1)}
	\end{split}
\end{align*}

where we use the fact that for two events $A,B$ we have $\mathbb P \big[ A \cap B \big] \geq 1-\mathbb P \big[ A^C\big]-\mathbb P \big[B^C \big]$ and Proposition \ref{prop:linear_implies_in_ball} and Proposition \ref{prop:axiom_implies_small_interesction} to conclude.

Now, we will use Lemma \ref{lem:implies_long_length} for $r=\kappa n$. For every $n \geq N_{\delta} / \kappa$  we get:

\begin{align*}
\begin{split}
	\mathbb P\Big[\dive(w^{p}_n,z^{p}_n,p, \delta)&> cn f\big(\frac{ \delta \kappa n}{4D}\big)\Big]  \geq  \mathbb P \left[\vert \mathcal A^{\delta \kappa n/2}_{p,z^{p}_n} \backslash \mathcal A^{\delta \kappa n/2}_{p,w^p_n} \vert-2 \geq cn \right]  \\
	&\geq 1-C_4e^{-n/C_4}-C_1e^{-\nu n/(2C_1)}.
\end{split}	
\end{align*}
Choosing a constant $C>0$ which also covers the cases $n < N_{\delta}/ \kappa$ and is such that $C \geq \max\{c^{-1}, 4D/\delta \kappa \}$ concludes the proof of Theorem \ref{thm:lower_bound}. \end{proof}


Finally, we can prove theorem \ref{thm:intro}.

\begin{proof}[Proof of Theorem \ref{thm:intro}]

We first look at the case when $(w_n)_n$ is a random walk. In this case, we need to make sure that the groups considered satisfy assumption \ref{assump:random_walk}, i.e. that they contain $f$-divergent elements for a suitable $f$. We then note in each case that the divergence of the whole group is equivalent to $nf(n)$ and hence Theorem \ref{thm:lower_bound} implies Theorem \ref{thm:intro}, in these cases.

\begin{itemize}
    \item In Example \ref{example:f_divgroups} (referring to \cite[Lemma 3.5]{GS21}) we show that relatively hyperbolic groups have exponential-divergent elements. In \cite[Theorem 1.3]{SistoRelativeHyp}, it is shown that a one-ended finitely presented relatively hyperbolic group have exponential divergence.
    \item  By \cite[Corollary 3]{HRSS} if $M$ is a 3-dimensional non-geometric graph manifold then $\pi_1(M)$ is a hierarchically hyperbolic group. Hence, by Lemma \ref{lem:HHS_divergent_element}, it contains a linear-divergent element. Now, by \cite{GerstenDivergence} graph manifolds have quadratic divergence.
    \item Mapping class groups are all examples of hierarchically hyperbolic groups \cite{HHS1}. Hence by Lemma \ref{lem:HHS_divergent_element}, these groups contain a linear-divergent element (when the corresponding surface is closed connected oriented and has genus at least 2). In all of these cases, the divergence is quadratic \cite{DuchinRafi}, \cite{BehrsotckMCG}. This is also implicit in \cite{BehrsotckMCG}.
\item Right-angled Artin groups are hierarchically hyperbolic groups with unbounded main curve graph when the defining graph is not a join \cite{HHS1}, hence by Lemma \ref{lem:HHS_divergent_element}, these groups contain a linear-divergent element. In the case where the defining graph is connected and not a join, the divergence is quadratic \cite{BehrstockCharney}. As a side note, there is also another reason for the existence of linear-divergent elements: By \cite[Theorem 5.2]{BehrstockCharney}, if $A_\Gamma$ is a right-angled Artin group whose defining graph $\Gamma$ is not a join, then $A_{\Gamma}$ contains a rank-one isometry for the action on the universal cover of its Salvetti complex $X_\Gamma$ which is CAT(0). Hence by Lemma \ref{lem:rank_one_divergent}, $A_{\Gamma}$ contains a linear-divergent element.
\end{itemize}

To deal with the case of tame Markov chains on groups from the first two items in Theorem \ref{thm:intro},  we need to show that Assumption \ref{assump:markov_chain} holds for these groups. We have just showed that these groups have relevant $f$-divergent elements. It remains to show that they satisfy the conclusion of \cite[Proposition 5.1]{GS21}. 

This is true as relatively hyperbolic groups satisfy \cite[Assumption 4.1]{GS21} by the fact that thay contain a superlinear-divergent element (which was called super-divergent element there) and graph manifold groups satisfy Assumption \cite[Assumption 4.3]{GS21}
(see the proof of \cite[Proposition 4.6]{GS21}).
\end{proof}

\bibliography{main}
\bibliographystyle{alpha}

\end{document}